\documentclass[11pt,reqno,a4paper,oneside]{amsart}

\usepackage[margin=1.75cm]{geometry}
\usepackage{latexsym}
\usepackage{amsmath}
\usepackage{amsthm}
\usepackage{amssymb}
\usepackage{graphicx}
\usepackage{caption}
\usepackage{subcaption}
\usepackage{amsfonts}
\usepackage{xcolor}
\usepackage{fancyhdr}
\usepackage{comment}
\usepackage[section]{placeins}
\usepackage{enumerate}
\usepackage{enumitem}

\usepackage{mathrsfs}                           

\newcommand{\eps}{\varepsilon}                
\newcommand{\bx}{\mathbf{x}}
\newcommand{\be}{\mathbf{e}}
\newcommand{\by}{\mathbf{y}}

\newcommand{\bzero}{\mathbf{0}}

\newcommand{\bs}{\mathbf{s}}
\newcommand{\mK}{\mathcal{K}}
\newcommand{\bk}{\mathbf{k}}
\newcommand{\bS}{\mathbf{S}}

\newcommand{\bbR}{\mathbb{R}}
\newcommand{\bbZ}{\mathbb{Z}}
\newcommand{\bi}{\mathbf{i}}
\newcommand{\bj}{\mathbf{j}}

\newcommand{\bm}{\mathbf{m}}
\newcommand{\bz}{\mathbf{z}}

\newcommand{\mG}{\mathcal{G}}

\newcommand{\mA}{\mathcal{A}}
\newcommand{\mAp}{\mathcal{A}_\alpha^{(p)}}
\newcommand{\mAn}{\mathcal{A}_\alpha^{(n)}}

\newcommand{\BE}{\begin{equation}}

\newcommand{\bes}{\begin{subequations}}
\newcommand{\EE}{\end{equation}}
\newcommand{\ees}{\end{subequations}}
\newcommand{\mO}{\mathcal{O}}

\newcommand{\BEU}{\begin{equation*}}
\newcommand{\EEU}{\end{equation*}}

\numberwithin{equation}{section}

\title[L\'{e}vy flight narrow capture problem]{Green's functions of the fractional Laplacian on a square - boundary considerations and applications to the L\'{e}vy flight narrow capture problem}
\author{J.~C.~Tzou$^\dagger$}

\thanks{$^\dagger$ School of Mathematical and Physical Sciences, Macquarie University, Sydney, NSW, Australia (tzou.justin@gmail.com)}
\thanks{\textbf{Funding:} JCT was supported by the Australian Research Council (DP220101808).}

\begin{document}
\begin{abstract}

On the unit square, we introduce a method for accurately computing source-neutral Green's functions of the fractional Laplacian operator with either periodic or homogeneous Neumann boundary conditions. This method involves analytically constructing the singular behavior of the Green's function in a neighborhood around the location of the singularity, and then formulating a ``smooth'' problem for the remainder term. This smooth problem can be solved for numerically using a basic finite difference scheme. This approach allows accurate extraction of the regular part of the Green's function (and its gradient, if so desired). This new tool enables quantification of properties and characteristics of the narrow capture problem, where a particle undergoing a L\'{e}vy flight of index $\alpha \in (0,1)$ searches for small target(s) of radius $\mO(\eps)$ for $0 < \eps \ll 1$ on a bounded two-dimensional domain. In particular, it allows us to show how boundary interactions and configuration of multiple targets impact expected search time. Furthermore, we are able to illustrate how a target can be ``shielded'' by obstacles, and how a L\'{e}vy flight search can be significantly superior in navigating these obstacles versus Brownian motion. All asymptotic predictions are confirmed by full numerical solutions.

\end{abstract}
\maketitle
\textbf{Keywords: fractional Laplacian, boundary conditions, Green's functions, L\'{e}vy flight, narrow capture problem, mean first passage time, splitting probabilities} 

\section{Introduction}

Since 2004 \cite{holcman2004escape}, the classic narrow capture/escape problem, which seeks the mean time (and related statistics) for a particle undergoing Brownian motion in a bounded (multi-dimensional) domain to first reach a set of small targets in the interior (narrow capture) or on the boundary (narrow escape) of the domain, has been used as a model for a wide array of applications. A small sampling includes the exit of a sodium ion through narrow valves on the cell membrane \cite{singer2006narrow1}; a diffusing intracellular molecule in search of a localized reaction site \cite{bressloff2022narrow}; animals foraging for food and water sources \cite{kurella2015asymptotic, mckenzie2009first}. We refer to the previous references, along with, e.g., \cite{holcman2015stochastic, holcman2014narrow}, for a more comprehensive review of applications that the narrow capture/escape problem has been used to model.

In the context of Brownian motion, with the Laplacian operator as the infinitesimal generator, this class of problems has been subject to analysis within various geometric settings. Asymptotic, microlocal, conformal mapping, and numerical methods have been used to consider flat two-dimensional geometries (e.g., \cite{pillay2010asymptotic, singer2006narrow1, singer2006narrow2}), cusp-like geometries \cite{guerrier2014brownian, holcman2011narrow, singer2006narrow3}, Riemannian surfaces \cite{coombs2009diffusion, nursultanov2023narrow, singer2006narrow3}, and three-dimensional geometries \cite{chen2011asymptotic, gomez2015asymptotic, cheviakov2010asymptotic}. 

Variations on the narrow capture/escape problem include cases in which targets are mobile \cite{tzou2015mean, lindsay2017optimization, iyaniwura2021simulation} or where they are partially absorbing, modeling scenarios in which an escape or reaction event occurs with probability less than one upon the particle encountering the target \cite{bressloff2022narrow, lindsay2015narrow, guerin2021universal, guerin2023imperfect, lindsay2017first}. Recently, stochastic resetting has been incorporated, where the diffusing particle's position resets to some fixed position in the domain according to a Poisson process \cite{bressloff2020search, bressloff2021asymptotic}. 

Related problems include computation of the principal eigenvalue of the Laplacian in various two-dimensional domains with Neumann \cite{kolokolnikov2005optimizing} and periodic boundary conditions \cite{paquin2020asymptotics}, computation of the variance of the first passage time (as opposed to its first moment) \cite{kurella2015asymptotic, lindsay2016hybrid}, full distributions of the first passage time (i.e., all moments) \cite{lindsay2016hybrid, grebenkov2019full, bressloff2020target, bressloff2021asymptotic}, and splitting probabilities \cite{delgado2015conditional, kurella2015asymptotic, cheviakov2010asymptotic}, i.e., the probability that the diffusing particle reaches a desired target \emph{before} hitting some other obstacle(s) in the search domain.

Many results for the above are given in terms of relevant Green's functions of the Laplacian operator on the domain in which the problem is posed. In some cases, such as the unit disk \cite{pillay2010asymptotic}, ellipse \cite{iyaniwura2021optimization}, or periodic Bravais lattices \cite{iyaniwura2021asymptotic}, explicit formulae for these Green's functions are known. In more general cases, they must be computed numerically via a scheme that yields sufficient accuracy for the particular application \cite{lindsay2016hybrid}.

While a fairly comprehensive suite of results for the narrow capture/escape framework has been compiled for the case when the searching particle undergoes Brownian motion, significantly less is known when the particle undergoes a L\'{e}vy flight \cite{chechkin2008introduction}. In contrast to the linear-in-time mean squared displacement $\langle |\Delta x|^2 \rangle \sim \Delta t$ of Brownian motion, the L\'{e}vy flight model of superdiffusion possesses an unbounded mean squared displacement. Instead, a L\'{e}vy flight of index $\alpha \in (0,1)$ in $d$ spatial dimensions is characterized by a jump distribution $p(\ell)$ with power-tail $p(\ell) \sim \ell^{-d - 2\alpha}$ resulting in a diverging variance and the superlinear scaling $\langle |\Delta x|^{2\delta} \rangle^{1/\delta} \sim \Delta t^{1/\alpha}$ for $\delta < \alpha$.

L\'{e}vy flight dynamics have been observed in animal foraging behavior, from insects to sharks to shearwaters \cite{viswanathan2011physics, humphries2012foraging}, prompting studies comparing the efficacy of random searches by a particle undergoing L\'{e}vy flight to one undergoing Brownian motion (see, e.g., \cite{lomholt2008levy, bartumeus2002optimizing}). L\'{e}vy-like characteristics have also been observed in bacterial motion \cite{ariel2015swarming,chechkin2008introduction}, and have served as a basis for numerical search algorithms \cite{yang2010eagle, heidari2017efficient}. Materials have also been developed through which the diffusive transport of light is governed by L\'{e}vy flights \cite{barthelemy2008levy}.

Here, we focus on a particular aspect of L\'{e}vy flight search on a bounded two-dimensional domain in which the target(s) is (are) small in comparison to the search domain. Analytic progress of this L\'{e}vy flight variant of the narrow capture problem has only begun in earnest recently: in \cite{tzou2024counterexample}, a two-term expansion was derived for the global mean first passage time (GMFPT) on the unit square and cube with periodic boundary conditions and a single target. The expansion for the unit square is given below in \eqref{ubar}, where the small parameter $0<\eps\ll 1$ is the target radius. In \cite{gomez2024first}, the analogous expansion was derived for multiple targets on a one-dimensional interval with periodic boundary conditions. In both results, the $\mO(1)$ correction term was expressed in terms of a certain Green's function that encodes information about the geometry of the search domain and/or the spatial distribution of targets in the domain. In \cite{chaubet2025geodesic}, a leading order estimate, similar to the leading order term of \cite{tzou2024counterexample} (see \eqref{ubar} below) and \cite{gomez2024first} was derived for Anosov Riemannian manifolds (without boundary). 

As we detail below, the primary new contributions in this paper are (1) the formulation of a Neumann-type boundary condition (in contrast to the periodic boundary conditions of \cite{tzou2024counterexample, gomez2024first}) in the special case of a square domain, and (2) the introduction of a method for computing Green's functions that does not rely on expansions in terms of explicitly known eigenfunctions. These new developments allow us to evaluate the aforementioned $\mO(1)$ correction terms and thereby illustrate boundary effects and shielding effects in the L\'{e}vy flight narrow capture problems on the unit square. We note that our formulation of the Neumann-type boundary condition extends also to rectangular domains, while our method for computing Green's functions with periodic boundary conditions generalizes to arbitrary periodic Bravais lattices.

In contrast to Brownian motion, the infinitesimal generator for a L\'{e}vy flight is nonlocal, which presents new difficulties in analysis as well as numerical computation over the classical narrow capture/escape problem. One consideration is that of formulation (and interpretation) of Neumann-type boundary conditions. In \cite{gomez2024first, tzou2024counterexample}, the boundaries were assumed to be periodic, while in \cite{chaubet2025geodesic}, the manifolds considered were without boundary. Most physical systems, however, are not accurately described by periodic boundary conditions -- indeed, narrow capture/escape problems are often formulated with homogeneous Neumann conditions assuming particle reflection at the boundary. We refer to \cite{parks2025nonlocal, dipierro2017nonlocal} for a discussion around challenges in prescribing Neumann-type boundary conditions in nonlocal frameworks, as well as an overview of those proposed. Another difficulty in the L\'{e}vy flight narrow capture/escape problem is that of computing the required Green's functions. In one dimension, \cite{gomez2024first} provides a rapidly converging infinite series in terms of eigenfunctions of the Laplacian with periodic boundary conditions, but no such series has been given in higher dimensions. In such cases, including those in which eigenfunctions are not explicitly known, an accurate numerical method must be developed for computing Green's functions.

In this paper, we address these two difficulties in the special case of the two-dimensional square (only slight modifications are needed to generalize to all rectangular domains). That is, we (1) address the formulation and physical interpretation of a Neumann-type boundary condition  discussed in \cite{andreu2010nonlocal,montefusco2012fractional,stinga2015fractional} for L\'{e}vy flights in a bounded domain, and (2) introduce a method for accurately computing Green's functions of the fractional Laplacian with both the Neumann-type as well as periodic boundary conditions. 

We apply these new developments within the context of the the L\'{e}vy flight narrow capture problem. For the calculation of the GMFPT, we highlight effects of reflecting boundary conditions and also of target configuration, while for the splitting probability problem, we illustrate ``shielding'' effects in which the target is surrounded by obstacles. Regarding the latter, we show how L\'{e}vy flight searches become less susceptible to such shielding effects as $\alpha$ decreases away from its Brownian limit of $1$. All of these effects are encoded in corrections terms involving certain Green's functions, which we derive and use our new method to compute. With respect to the boundary condition discourse, we emphasize here that our purpose is to illuminate a particularly simple physical interpretation of a Neumann-type boundary on a geometry that facilitates a convenient method for analysis and computation -- we again refer to \cite{parks2025nonlocal, dipierro2017nonlocal} for a thorough discussion on the various Neumann-type conditions proposed and their theoretical underpinnings.



Prior to outlining the sections of the paper, we give a brief overview of the main result of \cite{tzou2024counterexample} for the L\'{e}vy flight narrow capture problem with a single target of radius $0 < \eps \ll 1$ on the unit square $\Omega = [0,1]\times[0,1]$ with periodic boundary conditions. For a search conducted by a particle undergoing a L\'{e}vy flight of index $\alpha \in (0,1)$, with $\alpha \to 1^{-}$ being the Brownian limit, the GMFPT, $\bar{u}_\eps^{(p)}$, is given by
\BE \label{ubar}
	\bar{u}_\eps^{(p)} \sim \eps^{2\alpha-2}\frac{\Gamma(1-\alpha)}{4^\alpha\pi\Gamma(\alpha)} \chi_\alpha - R_\alpha^{(p)}(\bx_0;\bx_0) \,; \qquad \chi_\alpha \equiv \frac{\pi(1-\alpha)}{\sin[(1-\alpha)\pi)]} \,.
\EE
In \eqref{ubar}, the constant $\chi_\alpha$ was obtained through recasting a certain inner problem as an integral equation on the rescaled domain of the target; for the disk-shaped target considered in \cite{tzou2024counterexample}, the integral equation yielded an analytic solution for $\chi_\alpha$ given by \cite{kahane1981solution}. The $\mO(1)$ correction term in \eqref{ubar}, $R_\alpha^{(p)}(\bx_0;\bx_0)$, is the regular part of a certain Green's function $G_\alpha^{(p)}(\bx;\bx_0)$ periodic on $\Omega$ satisfying 
\bes \label{Gperiodic}
\BE \label{Gperiodiceq}
	\mathcal{A}_\alpha^{(p)} G_\alpha^{(p)} = - 1 + \delta(\bx-\bx_0) \,, \quad \bx \in \Omega \setminus \lbrace\bx_0\rbrace \,; \qquad \int_{\Omega} \! G_\alpha^{(p)}(\bx;\bx_0)\, d\bx = 0 \,,
\EE
\BE \label{Gperiodicloc}
	G_\alpha^{(p)} (\bx;\bx_0) \sim - \frac{c_\alpha}{|\bx-\bx_0|^{2-2\alpha}} + R_\alpha^{(p)}(\bx;\bx_0) + \mO(|\bx-\bx_0|) \enspace \mbox{as} \enspace \bx \to \bx_0 \,; \qquad c_\alpha \equiv \frac{\Gamma(1-\alpha)}{4^\alpha\pi\Gamma(\alpha)} \,.
\EE
\ees

Note that the required periodic boundary conditions on $\partial\Omega$ are encoded in the operator $\mathcal{A}_\alpha^{(p)}$ itself, derived in \cite{tzou2024counterexample} through a discrete-continuum probabilistic argument and given by
\bes
\BE \label{Ap}
\mathcal{A}_\alpha^{(p)} f(\bx) \equiv C_\alpha  \int_\Omega \!  [f(\by)-f(\bx)] \sum_{\bm \in \mathbb{Z}^2} \frac{1}{|T_\bm^{(p)}(\by) - \bx|^{2+2\alpha}} \, d\by \,, \quad \bx = (x_1,x_2) \in \Omega;
\EE
\BE \label{Tm}
	T_\bm^{(p)}(\by)  = \by + \bm \,, \quad \bm \in \bbZ^2 \,,
\EE
for $f(\bx)$ periodic on $\Omega$. In \eqref{Ap}, the constant $C_\alpha$ is defined by
\BE \label{Calpha}
	C_\alpha = \frac{4^\alpha\Gamma(1+\alpha)}{\pi|\Gamma(-\alpha)|} \,,
\EE
\ees
which ensures that $A_\alpha{(p)}\phi_i^{(p)} = -|\lambda_i^{(p)}|^{\alpha}\phi_i^{(p)}$, where $(\lambda_i^{(p)}, \phi_i^{(p)}(\bx))$ are eigenpair of the regular Laplacian.

Regarding \eqref{ubar}, the leading order term was confirmed in \cite{tzou2024counterexample} via both Monte Carlo simulations as well as a finite difference solution of the narrow escape elliptic problem $\mathcal{A}_\alpha^{(p)} u_\eps^{(p)} = -1$ with $u = 0$ when $|\bx-\bx_0| < \eps$. However, in the absence of a method for accurately computing the regular part $R_\alpha^{(p)}(\bx_0;\bx_0)$, the $\mO(1)$ term was left unverified. Furthermore, the periodic boundary conditions assumed in \cite{tzou2024counterexample}, which are encoded in $\mathcal{A}_\alpha^{(p)}$ and which assume that a L\'{e}vy flight particle exiting $\Omega$ on one edge via a long jump returns through the opposite edge, do not lend well to the modeling of physical systems. 

We address these two shortcomings in \S\S \ref{sec:reflect} and \ref{sec:Greens}. In \S \ref{sec:reflect}, we construct an operator $\mathcal{A}_\alpha^{(n)}$ on $\Omega = [0,1]\times[0,1]$ that models specular reflection at $\partial\Omega$ of particles that hit $\partial\Omega$ mid-jump and would have otherwise exited the domain. Through a heuristic argument, we show that this operator is in fact equivalent to the spectral representation of the fractional Laplacian on $\Omega$ with restriction to the set of Laplacian eigenfunctions $\phi_i^{(n)}$, $i = 0, 1, \ldots$, satisfying $\partial_\nu \phi_i^{(n)} = 0$ on $\partial \Omega$, where $\partial_\nu$ denotes the normal derivative. That is, for these Laplacian eigenfunctions on $\Omega$ corresponding to eigenvalues $\lambda_i^{(n)} \leq 0$, $i = 0, 1, \ldots$, respectively, we argue that $\mathcal{A}_\alpha^{(n)}\phi_i^{(n)} = -|\lambda_i^{(n)}|^{\alpha} \phi_i^{(n)}$. We note that the specification of the usual $\partial_\nu \phi_i^{(n)} = 0$ condition on $\partial\Omega$ through the spectral representation of the fractional Laplacian was discussed in \cite{andreu2010nonlocal,montefusco2012fractional,stinga2015fractional}.

In \S \ref{sec:Greens}, we introduce a method for accurate computation of the periodic Green's function $G_\alpha^{(p)}$ of \eqref{Gperiodic} as well as the analogous Neumann Green's function $G_\alpha^{(n)}$ for the $\mathcal{A}_\alpha^{(n)}$ operator with $\partial_\nu G_\alpha^{(n)} = 0$ on $\partial\Omega$. By constructing the singular part of the Green's functions analytically, and formulating a smooth numerical problem for the remainder term, this method yields accurate values for the regular parts $R_\alpha^{(p)}(\bx;\bx_0)$ for the periodic Green's function and its Neumann counterpart $R_\alpha^{(n)}(\bx;\bx_0)$. We note that this method does not invoke the spectral interpretation of $\mathcal{A}_\alpha^{(p)}$ and  $\mathcal{A}_\alpha^{(n)}$, and thus does not rely on possession of the Laplacian eigenfunctions on $\Omega$.

\begin{figure}
	\centering
	\begin{subfigure}{0.49\textwidth}
		\centering
		\includegraphics[width=1.1\textwidth]{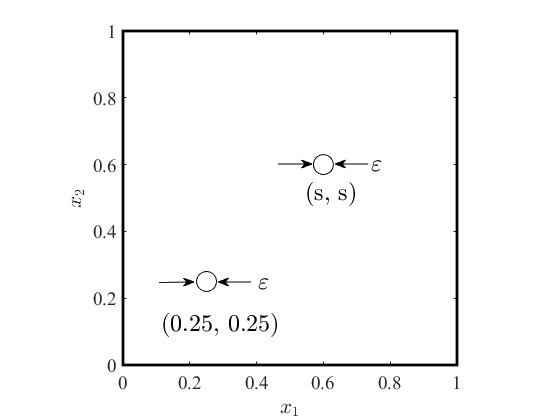}
		\caption{targets at $(0.25, 0.25)$ and $(s,s)$ - periodic BCs}
		\label{fig:schematic_twotraps}
	\end{subfigure}
	\begin{subfigure}{0.49\textwidth}
		\centering
		\includegraphics[width=1.1\textwidth]{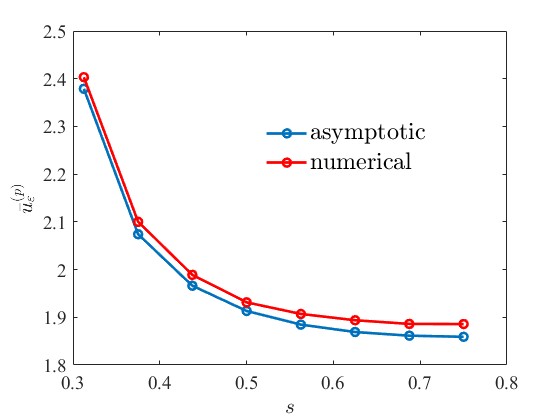}
		\caption{GMFPT versus $s$}
		\label{fig:mfpt_vs_second_trap_location_asympt_vs_numerical_periodic}
	\end{subfigure}
	\caption{(a) In $\Omega$ with periodic boundary conditions, we center one circular target at $(0.25, 0.25)$ and another at $(s,s)$, where $0.25<s<0.75$. Both targets have radius $0<\varepsilon \ll 1$. (b) The global mean first passage time versus $s$ for $\alpha = 0.6$ and $\varepsilon = 0.03$. The red curve is generated from a finite difference solution for $u_\eps^{(p)}$ satisfying \eqref{narrowescapeequationmult} with $\mA$ replaced by $\mAp$, while the blue curve is obtained through an asymptotic analysis along with the algorithm for accurate computation of $G_\alpha^{(p)}$. }
	\label{fig:bcs}
\end{figure}

In \S \ref{sec:narrow}, we leverage the Green's function computations of \S \ref{sec:Greens} to compute the $\mO(1)$ correction terms of the GMFPT. For the periodic domain, we focus on the case of multiple targets, deriving the $\mO(1)$ correction term in terms of the Green's function $G_\alpha^{(p)}$ of \eqref{Gperiodic} to predict the effect of target configuration on GMFPT. The schematic of this problem is shown in Fig. \ref{fig:schematic_twotraps}, where two targets of radius $\eps$ are centered at $(0.25, 0.25)$ and $(s,s)$, respectively. In Fig. \ref{fig:mfpt_vs_second_trap_location_asympt_vs_numerical_periodic}, we plot the GMFPT as $s$ is varied, showing that it reaches a minimum when $s = 0.75$. We note the close agreement between the asymptotic prediction and the numerical values.

\begin{figure}
	\centering
	\begin{subfigure}{0.49\textwidth}
		\centering
		\includegraphics[width=1.1\textwidth]{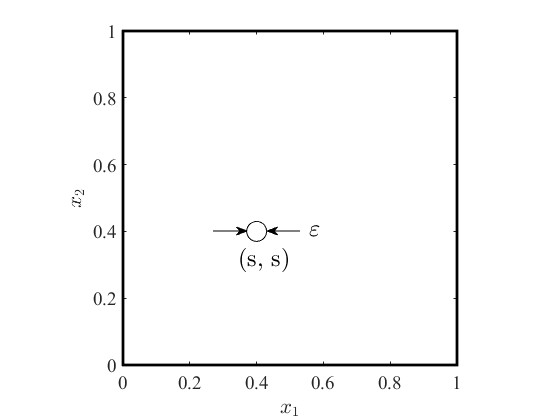}
		\caption{target centered at $(s,s)$ - reflective BCs}
		\label{fig:schematic_onetrap}
	\end{subfigure}
	\begin{subfigure}{0.49\textwidth}
		\centering
		\includegraphics[width=1.1\textwidth]{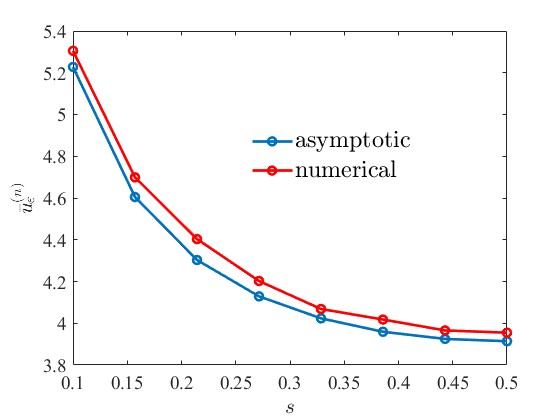}
		\caption{GMFPT versus $s$}
		\label{fig:mfpt_vs_trap_location_asympt_vs_numerical_neumann}
	\end{subfigure}
	\caption{(a) In $\Omega$ with reflective boundary conditions, we center a single circular target of radius $0<\varepsilon \ll 1$ at $(s,s)$, with $s \in (0,0.5)$. (b) The global mean first passage time versus $s$ for $\alpha = 0.6$ and $\varepsilon = 0.03$. The red curve is generated from a finite difference solution for $u_\eps^{(n)}$ satisfying \eqref{narrowescapeequationmult} with $\mA$ replaced by $\mAn$, while the blue curve is obtained through an asymptotic analysis along with the algorithm for accurate computation of $G_\alpha^{(n)}$.}
	\label{fig:bcs}
\end{figure}

Due to the symmetry of a periodic domain, the GMFPT is independent of target location in the single-target problem. This symmetry however, is broken in the case of reflective boundary conditions. For the case of specular reflection at the boundary, we use our computation of the Green's function $G_\alpha^{(n)}$ to predict the how GMFPT changes as the location of a single target is varied in $\Omega$. The schematic of this problem is shown in Fig. \ref{fig:schematic_onetrap}, where one target of radius $\eps$ is centered at $(s, s)$. In Fig. \ref{fig:mfpt_vs_trap_location_asympt_vs_numerical_neumann}, we plot the GMFPT as $s$ is varied, showing that it reaches a minimum when $s = 0.5$. We note the close agreement between the asymptotic prediction and the numerical values. This variation in the GMFPT as the target location is varied is absent in the case of periodic boundary conditions, and is due solely to the effects of the reflective boundary.

\begin{figure}
	\centering
	\begin{subfigure}{0.49\textwidth}
		\centering
		\includegraphics[width=1.1\textwidth]{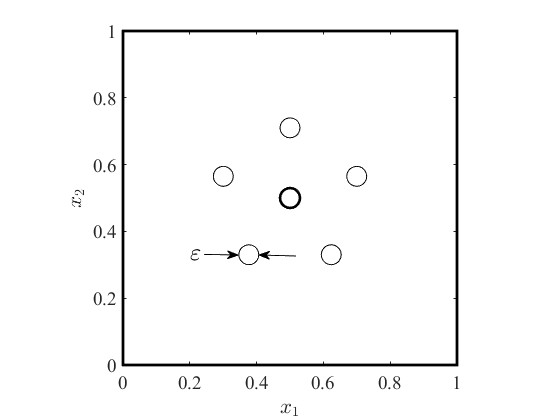}
		\caption{target centered at $(s,s)$ - reflective BCs}
		\label{fig:schematic_splitting}
	\end{subfigure}
	\begin{subfigure}{0.49\textwidth}
		\centering
		\includegraphics[width=1.1\textwidth]{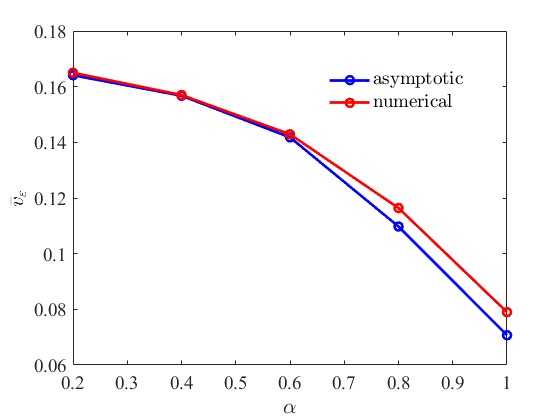}
		\caption{average splitting probability versus $\alpha$}
		\label{fig:split_prob_vs_alpha}
	\end{subfigure}
	\caption{(a) In $\Omega$ with reflective boundary conditions, the desired target at the center (heavy line) is ``shielded'' by five obstacle targets. The average splitting probability, $\bar{v}_\eps^{(n)}$, is the probability of reaching the desired target \emph{before} hitting any of the obstacle targets, averaged over all starting locations in $\Omega$. (b) Plot of $\bar{v}_\eps^{(n)}$ versus $\alpha$, where $\alpha = 1$ is the Brownian limit. The red curve is generated from a finite difference solution for $v_\eps^{(n)}$ satisfying \eqref{splitting}, while the blue curve is obtained through an asymptotic analysis along with the algorithm for accurate computation of $G_\alpha^{(n)}$. As expected, L\'{e}vy flights with smaller index $\alpha$, which experience more long jumps, are less susceptible to the shielding effect.}
	\label{fig:bcs}
\end{figure}

In \S \ref{sec:split}, on the unit square with reflective boundaries, we consider the splitting probability of reaching a desired target \emph{before} hitting one of the obstacle targets that surround it. The schematic of this problem is depicted in Fig. \ref{fig:schematic_splitting}, where the desired target (heavy line) is ``shielded'' by five other targets (light line). In particular, we compare this splitting probability across various  L\'{e}vy flight indices $\alpha$ ranging from $0.2$ to the Brownian limit of $\alpha = 1$. The results are shown in Fig. \ref{fig:split_prob_vs_alpha}.
 Through an asymptotic analysis of a certain elliptic problem, we derive an expansion for the splitting probability in terms of quantities associated with the Neumann Green's function $G_\alpha^{(n)}$. We demonstrate that as $\alpha$ increases, the average probability of reaching the desired target prior to one of the surrounding targets decreases, illustrating the (perhaps expected) phenomenon that a L\'{e}vy flights with its long jumps are less susceptible to ``shielding'' effects than Brownian motion.  In \S \ref{sec:disc}, we discuss other possible applications of our work as well as how it can be extended. We also comment on the limitations of our work, and conclude with a brief discussion of some related open problems.
 
\section{The fractional Laplacian for reflective boundary conditions on $\partial\Omega$} \label{sec:reflect}

\begin{figure}
\centering
\begin{subfigure}{0.49\textwidth}
	\centering
	\includegraphics[width=1.1\textwidth]{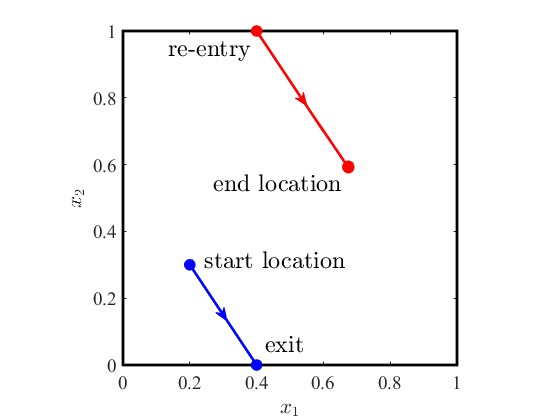}
	\caption{periodic boundary condition}
	\label{fig:periodicbc}
\end{subfigure}
\begin{subfigure}{0.49\textwidth}
	\centering
	\includegraphics[width=1.1\textwidth]{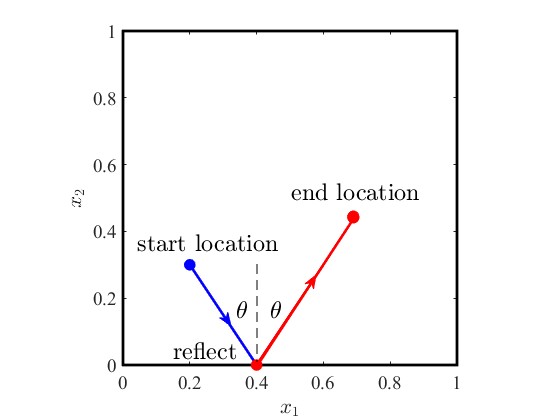}
	\caption{specular reflection at boundary}
	\label{fig:neumannbc}
\end{subfigure}
\caption{Periodic boundary (a) versus specular reflection at the boundary (b). In (a), a particle whose trajectory causes it to exit the domain through the bottom boundary simply re-enters at the same angle from the top boundary. In (b), the same trajectory is reflected off the boundary at the same angle as the incident trajectory. The total distance traveled by the particle is the same in both cases; the distribution of this distance follows a power-law distribution.}
\label{fig:bcs}
\end{figure}

In this section we derive an expression for the operator $\mAn$ describing a L\'{e}vy flight process with index $\alpha \in (0,1)$ in the domain $\Omega = [0,1]\times[0,1]$, where we assume specular reflection of particles on $\partial \Omega$. That is, we assume that a particle whose jump would take it outside of $\Omega$ instead reflects off $\partial\Omega$ in a specular manner and remains inside $\Omega$ (see Fig. \ref{fig:neumannbc}).  This is in contrast to the periodic boundary conditions depicted in Fig. \ref{fig:periodicbc}, where a particle exiting the domain simply re-enters through the opposite edge. 

For $\mAn$ derived under this assumption of specular reflection on $\partial \Omega$, we argue that $\mAn\phi_i^{(n)} = -|\lambda_i^{(n)}|^{\alpha}\phi_i^{(n)}$, $i = 0, 1, \ldots$. Here, $\phi_i^{(n)}$ is the $i$-th eigenfunction of the Laplacian operator on $\Omega$ with $\partial_\nu \phi_i^{(n)} = 0$ on $\partial\Omega$ while $\lambda_i^{(n)} \leq 0$ is the corresponding eigenvalue. To see this, let $\Phi_i^{(n)}(\bx)$, $\bx \in \mathbb{R}^2$, be the even extension of $\phi_i^{(n)}$ onto $\mathbb{R}^2$. That is, suppose $\by = (y_1,y_2) \in \Omega$, and let 
\BE \label{Y}
T_\bm^{(n)}(\by) = (m_1 + \!\!\!\!\!\!\ \mod(m_1,2) + (-1)^{m_1} y_1, m_2 + \!\!\!\!\!  \mod(m_2,2) + (-1)^{m_2} y_2) \,; \qquad \bm = (m_1,m_2) \in \mathbb{Z}^2 \,.
\EE
Then  $\Phi_i^{(n)}(T_\bm^{(n)}(\by)) = \phi_i^{(n)}(\by)$. In \eqref{Y}, $T_\bm^{(n)}(\by)$ are simply image points of $\mathbf{y}$; the first component is obtained by successively reflecting $\mathbf{y}$ a number $|m_1|$ times across the lines $1, 2, \ldots, m_1$ when $m_1 > 0$ and across the lines $0, -1, \ldots, m_1 + 1$ when $m_1 < 0$, and similarly for the second component. The point $T_\bm^{(n)}(\by)$ lies in the unit square whose lower left vertex is at the point $\bm$; we denote this unit square $\Omega_\bm$ so that $\Omega_\bzero \equiv \Omega$. Note that $T_\bm^{(n)}(\by)$ is a one-to-one map from $\Omega$ to $\Omega_\bm$. 

For the fractional Laplacian $-(-\Delta)^{\alpha}$ defined on $\mathbb{R}^2$, we have that $-(-\Delta)^{\alpha} \Phi_i^{(n)}(\bx)$ is given by
\BE \label{fraclap}
	-(-\Delta)^{\alpha} \Phi_i^{(n)}(\bx) \equiv C_\alpha \int_{\mathbb{R}^2} \! \frac{\Phi_i^{(n)}(\by) - \Phi_i^{(n)}(\bx)}{|\by-\bx|^{2+2\alpha}} \, d\by = -|\lambda_i^{(n)}|^{\alpha}\Phi_i^{(n)}(\bx) \,; \qquad \bx \in \mathbb{R}^2 \,,
\EE
where the integral near $\by = \bx$ is taken in the principal value sense. Since $\mathbb{R}^2$ can be tiled as $\bigcup\limits_{\bm \in \mathbb{Z}^2} \Omega_\bm$, we may rewrite the right-hand side of \eqref{fraclap} as 
\BE \label{fraclapsum}
-(-\Delta)^{\alpha} \Phi_i^{(n)}(\bx) = C_\alpha \sum_{\bm\in\mathbb{Z}^2} \int_{\Omega_\bm} \! \frac{\Phi_i^{(n)}(\by) - \Phi_i^{(n)}(\bx)}{|\by-\bx|^{2+2\alpha}} \, d\by = -|\lambda_i^{(n)}|^{\alpha}\Phi_i^{(n)}(\bx) \,; \qquad \bx \in \mathbb{R}^2 \,.
\EE
Since each $\by \in \Omega_\bm$ has one corresponding point in $\Omega$, we make the change of variable $\by \to T_{\mathbf{m}}(\by)$ in \eqref{fraclapsum} for each integration over $\Omega_\bm$
\BE \label{fraclapsum2}
 -(-\Delta)^{\alpha} \Phi_i^{(n)}(\bx) = C_\alpha \sum_{\bm\in\mathbb{Z}^2}\int_{\Omega} \! \frac{\Phi_i^{(n)}(T_{\mathbf{m}}(\by)) - \Phi_i^{(n)}(\bx)}{|T_\bm^{(n)}(\by)-\bx|^{2+2\alpha}} \, d\by = -|\lambda_i^{(n)}|^{\alpha}\Phi_i^{(n)}(\bx) \,; \qquad \bx \in \mathbb{R}^2 \,,
 \EE
where $T_{\mathbf{m}}(\by)$ is given in \eqref{Y}. Using that  $\Phi_i^{(n)}(T_{\mathbf{m}}(\by)) = \phi_i^{(n)}(\by)$, and restricting $\bx$ to $\Omega$, we obtain from \eqref{fraclapsum2}
\BE \label{fraclapsum3}
-(-\Delta)^{\alpha}\phi_i^{(n)}(\bx) = C_\alpha \sum_{\bm\in\mathbb{Z}^2}\int_{\Omega} \! \frac{\phi_i^{(n)}(\by) - \phi_i^{(n)}(\bx)}{|T_\bm^{(n)}(\by)-\bx|^{2+2\alpha}} \, d\by = -|\lambda_i^{(n)}|^{\alpha}\phi_i^{(n)}(\bx) \,; \qquad \bx \in \Omega \,.
\EE
Interchanging the sum and integral in \eqref{fraclapsum3}, we obtain
\BE \label{fraclapsum4}
\mAn\phi_i^{(n)}(\bx) \equiv C_\alpha \int_{\Omega} \! \left[\phi_i^{(n)}(\by) - \phi_i^{(n)}(\bx)\right] \sum_{\bm\in\mathbb{Z}^2}  \frac{1}{|T_\bm^{(n)}(\by)-\bx|^{2+2\alpha}} \, d\by = -|\lambda_i^{(n)}|^{\alpha}\phi_i^{(n)}(\bx) \,; \qquad \bx \in \Omega \,.
\EE
Thus, with $T_{\mathbf{m}}(\by)$ defined in \eqref{Y}, the action of the operator $\mAn$, defined in \eqref{fraclapsum4} on Neumann eigenfunctions of the Laplacian on $\Omega$ is equivalent to that of the spectral definition of the fractional Laplacian on $\Omega$. Like its periodic counterpart, $\mAp$, it inherits the property of self-adjointness from the fractional Laplacian $-(-\Delta)^\alpha$ of \eqref{fraclap}, the latter of which is due the symbol of  $-(-\Delta)^\alpha$  being real-valued.

Now consider a function $u_\eps^{(n)}(\bx)$, $\bx \in \Omega$, that represents the mean first passage time (MFPT) of a L\'{e}vy flight starting from location $\bx$ to a target of radius $\eps$ centered at $\bx_0 \in \Omega$. Also, let $\partial_\nu u_\eps^{(n)} = 0$ on $\partial\Omega$. We claim that, for such a L\'{e}vy flight with index $\alpha \in (0,1)$, $u_\eps^{(n)}(\bx)$ satisfies 
\BE \label{narrowescapeequation}
	\mAn u_\eps^{(n)}(\bx) = -1 \,, \quad  \bx \in \Omega \setminus B_\eps(\bx_0) \,; \qquad  u_\eps^{(n)} = 0 \enspace \mbox{for} \enspace \bx \in  B_\eps(\bx_0) \,; \quad \bx_0 \in \Omega \,,
\EE
where $\mAn$ is the operator defined in \eqref{fraclapsum4}, and $B_\eps(\bx_0)$ denotes a disk of radius $\eps$ centered at $\bx_0$. Furthermore, we claim that \eqref{narrowescapeequation} describes specular reflection of the L\'{e}vy flight particle on $\partial \Omega$ as depicted in Fig. \ref{fig:neumannbc}.

To show this, let us first consider the even extension onto $\mathbb{R}^2$ of $u_\eps^{(n)}(\bx)$, which we denote $U_\eps^{(n)}(\bx)$. That is, for $\bx \in \Omega$ and $\bm \in \bbZ^2$, let $U_\eps^{(n)}$ be such that $U_\eps^{(n)}(T_\bm(\bx)) = u_\eps^{(n)}(\bx)$, where $T_\bm(\bx) \in \bbR^2$ is defined in \eqref{Y}. Note that $U_\eps^{(n)}(\bx)$ represents the MFPT from $\bx \in \bbR^2$ to the set of targets $\bigcup_{\bm \in \bbZ^2} B_\eps(T_\bm^{(n)}(\bx_0))$. 

We formulate the problem for $U_\eps^{(n)}(\bx)$ using the discrete-continuum approach employed in \cite{tzou2024counterexample} and \cite{valdinoci2009long} as follows. First, we discretize $\bbR^2$ into the set of discrete points $h\bi$ with $\bi \in \bbZ$, $h = 1/N$, $N \in \bbZ^+$, and $N \gg 1$. Then, by conditioning on the first jump, we have that the MFPT starting from point $h\bi$ is the weighted average of the MFPT's from all the points to which the particle can jump, plus the $\Delta t$ time it takes to make the jump:
\bes \label{constituentwhole}
\begin{equation} \label{constituent}
	U_\eps^{(n)}(h\bi) = \sum_{\bj \in \bbZ^2} w(h\bi, h\bj) U_\eps^{(n)}(h\bj) + \Delta t \,, \quad h\bi \in \bbR^2 \setminus \bigcup_{\bm \in \bbZ^2} B_\eps(T_\bm^{(n)}(\bx_0)) \,,
\end{equation}
with the exterior condition
\begin{equation} \label{constituentexterior}
	U_\eps^{(n)}(h\bi) = 0 \enspace \mbox{when} \enspace h\bi \in \bigcup_{\bm \in \bbZ^2} B_\eps(T_\bm^{(n)}(\bx_0)) \,.
\end{equation}
\ees
In \eqref{constituent}, the weight $w(h\bi, h\bj)$ is the probability of jumping from $h\bi$ to $h\bj$. Here, we take the approach of  Valdinoci in \cite{valdinoci2009long} and use the discrete power law distribution
\bes \label{powerlaw}
\BE \label{powerlawprob}
	w(h\bi,h\bj) = \begin{cases} 
		0 & \bi = \bj  \,,\\
		N_\alpha|\bi-\bj|^{-2-2\alpha} & \bi \neq \bj \,,
	\end{cases}
\EE
where $N_\alpha$ in \eqref{powerlawprob} is a normalization constant given by
\BE \label{powerlawprob2}
N_\alpha = \frac{1}{\sum\limits_{\bi = \bbZ^2, \bi \neq \bzero}|\bi|^{-2-2\alpha} } \,.
\EE
\ees
Since $\sum_{\bj \in \bbZ^2} w(h\bi,h\bj) = 1$ by \eqref{powerlaw}, we may rewrite \eqref{constituent} as
\begin{equation} \label{constituent2}
  \frac{N_\alpha}{\Delta t} \sum_{\bj \in \bbZ^2} \frac{U_\eps^{(n)}(h\bj) - U_\eps^{(n)}(h\bi)}{|\bi-\bj|^{2+2\alpha}} = -1 \,.
\end{equation}
Using the formal scaling law of \cite{valdinoci2009long}, $\Delta t = D_\alpha h^{2\alpha}$ for some constant $D_\alpha$, we obtain the Riemann sum on the left-hand side
\begin{equation} \label{constituent3}
	\frac{N_\alpha}{D_\alpha} h^{2}\sum_{\bj \in \bbZ^2} \frac{U_\eps^{(n)}(h\bj) - U_\eps^{(n)}(h\bi)}{|h\bi-h\bj|^{2+2\alpha}} = -1 \,.
\end{equation}
Letting $h\bi \equiv \bx$ and $h\bj \equiv \by$, \eqref{constituent3} together with \eqref{constituentexterior} yield the exterior problem
\bes
\BE \label{narrowcaptureR2}
	-(-\Delta)^{\alpha} U_\eps^{(n)}(\bx) = -1 \,, \quad \bx \in \bbR^2 \setminus \bigcup_{\bm \in \bbZ^2} B_\eps(T_\bm^{(n)}(\bx_0)) \,,
\EE
\BE
U_\eps^{(n)}(\bx) = 0 \enspace \mbox{when} \enspace \bx \in \bigcup_{\bm \in \bbZ^2} B_\eps(T_\bm^{(n)}(\bx_0)) \,,
\EE
\ees
where in \eqref{constituent3}, we have chosen the constant $D_\alpha$ so that
\[
\frac{N_\alpha}{D_\alpha} = \frac{4^\alpha\Gamma(1+\alpha)}{\pi|\Gamma(-\alpha)|} \,.
\]
Finally, restricting $\bx$ to $\Omega$ in \eqref{narrowcaptureR2}, and using the reflection symmetry $U_\eps^{(n)}(T_\bm^{(n)}(\bx)) = u_\eps^{(n)}(\bx)$, we have that $u_\eps^{(n)}(\bx)$ must satisfy \eqref{narrowescapeequation} with $\mAn$ defined in \eqref{fraclapsum4}.

\begin{figure}
	\centering
		\includegraphics[width=0.5\textwidth]{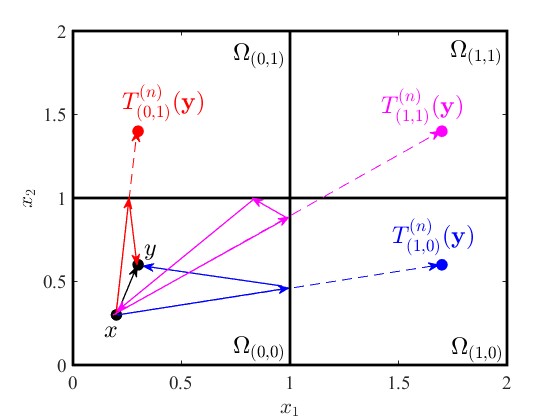}
	\caption{Four different paths from $\bx$ to $\by$ in $\Omega = [0,1]\times[0,1]$ with reflective boundaries. The direct path from $\bx$ to $\by$ of length $|\by-\bx|$ has probability $\sim|\by-\bx|^{-2-2\alpha}$. The other paths involve one or more reflections off of $\partial\Omega$, and have probability $\sim|T_{\bm}(\by)~-~\bx|^{-2-2\alpha}$ for $\bm \in \bbZ^2$.}
	\label{fig:reflection_fig}
\end{figure}

We now give a brief argument for why $u_\eps^{(n)}(\bx)$ satisfying \eqref{narrowescapeequation} is equivalent to the case of specular reflection on $\partial\Omega$. To see this, note that in the same way that the $|\bi-\bj|^{-2-2\alpha}$ term in \eqref{constituent3} is directly related to the probability of the particle jumping from $h\bi$ to $h\bj$, the infinite sum in $\sum_{\bm \in \bbZ^2}~|T_\bm^{(n)}(\by)~-~\bx|^{-2-2\alpha}$ arises from the infinite number of paths from $\bx$ to $\by$ and their associated probabilities. For example, the $\bm = \bzero$ term is simply the straight-line path from $\bx$ to $\by$ of length $|\by-\bx|$; the $\bm = (1,0)$ term represents the unique path from $\bx$ to $\by$ that involves a single specular reflection off of the $x_1 = 1$ boundary; the $\bm = (1,1)$ represents the unique path from $\bx$ to $\by$ that involves a single specular reflection off of each of the $x_1 = 1$ and $x_2 = 1$ boundaries, and so on. Sample paths from $\bx$ to $\by$, along with their trajectories, are shown in Fig. \ref{fig:reflection_fig}.


\section{Computation of the Green's functions} \label{sec:Greens}

In this section, we introduce a scheme for accurate computation of the regular part of the periodic Green's function $G_\alpha^{(p)}$ satisfying \eqref{Gperiodic} as well as a Neumann Green's function $G_\alpha^{(n)}$ satisfying
\bes \label{Gneumann}
\BE \label{Gneumanneq}
	\mAn G_\alpha^{(n)}(\bx;\bx_0) = -1 + \delta(\bx-\bx_0) \,, \qquad \bx \in \Omega \setminus \lbrace \bx_0 \rbrace \,; \qquad \int_{\Omega} \! G_\alpha^{(n)}(\bx;\bx_0)\, d\bx = 0 \,,
\EE
\BE \label{Gneumannloc}
G_\alpha^{(n)} (\bx;\bx_0) \sim - \frac{c_\alpha}{|\bx-\bx_0|^{2-2\alpha}} + R_\alpha^{(n)}(\bx;\bx_0) + \mO(|\bx-\bx_0|) \enspace \mbox{as} \enspace \bx \to \bx_0 \,,
\EE
\ees
for $\mAn$ defined in \eqref{fraclapsum4} and constant $c_\alpha$ defined in \eqref{Gperiodicloc}. Since the approach for both is similar, we drop the superscripts $^{(p)}$ and $^{(n)}$ for the remainder of this section and simply let $G_\alpha$ satisfy
\bes \label{G}
\BE \label{Geq}
\mathcal{A}_\alpha G_\alpha(\bx;\bx_0) = -1 + \delta(\bx-\bx_0) \,, \qquad \bx \in \Omega \setminus \lbrace \bx_0 \rbrace \,; \qquad \int_{\Omega} \! G_\alpha(\bx;\bx_0)\, d\bx = 0 \,,
\EE
\BE \label{Gloc}
G_\alpha (\bx;\bx_0) \sim - \frac{c_\alpha}{|\bx-\bx_0|^{2-2\alpha}} + R_\alpha(\bx;\bx_0) + \mO(|\bx-\bx_0|) \enspace \mbox{as} \enspace \bx \to \bx_0 \,.
\EE
\ees
We note that the right-hand side of \eqref{Geq} is orthogonal to the (co)kernel of the formally self-adjoint operator $\mA_\alpha$ (i.e., the set of all constant functions). The equation for $G_\alpha$ is thus consistent; the integral equation in \eqref{Geq} enforces uniqueness. 

We begin by decomposing the Green's function $G_\alpha(\bx;\bx_0)$ as
\bes
\begin{equation} \label{Gdecompose}
	G_\alpha(\bx;\bx_0) = \chi(\bx-\bx_0)u_0(\bx-\bx_0)  + \tilde{R}_\alpha(\bx;\bx_0) \,; \qquad u_0(\bx) \equiv -\frac{c_\alpha}{|\bx|^{2-2\alpha}} \,,
\end{equation}
where $\tilde{R}_\alpha(\bx;\bx_0)$ satisfies the same boundary conditions as $G_\alpha$, and is everywhere finite and infinitely smooth on $\Omega$ (see \cite{tzou2024counterexample}), while $\chi(\bx)$ is an infinitely smooth radially symmetric cut-off function, monotonic in $|\bx|$, centered at $\bx = \bzero$ such that, for $r_0 < r_1$,
\begin{equation} \label{cutoff}
   \chi(\bx) =  \begin{cases} 
		1 & 0\leq |\bx| < r_0 \\
	    0 & r_1 \leq |\bx| 
	\end{cases} \,,
\end{equation}
\ees
and $r_1$ is chosen so that the support of $\chi(\bx-\bx_0)$ lies entirely in $\Omega$. In \eqref{Gdecompose}, $u_0$ is the free-space Green's function of the fractional Laplacian satisfying
\BE \label{Gfree}
-(-\Delta)^{\alpha} u_0(\bx-\bx_0) = \delta(\bx-\bx_0) \,.
\EE
This can be easily seen by noting that $(-\Delta)^{\alpha}$ has the Fourier symbol $|\bk|^{2\alpha}$ along with the distributional Fourier transform pair $\mathcal{F}\left\lbrace c_\alpha|\bx|^{2\alpha-2}\right\rbrace(\bk) = |\bk|^{-2\alpha}$ in $\mathbb{R}^2$ for $0<\alpha<1$, where $c_\alpha$ is defined in \eqref{ubar}.

We will formulate a ``smooth'' equation for $\tilde{R}_\alpha(\bx;\bx_0)$ that can be solved for numerically using standard finite difference methods. Note that, comparing \eqref{Gloc} with \eqref{Gdecompose} and \eqref{cutoff}, we have $\tilde{R}_\alpha(\bx;\bx_0) = R_\alpha(\bx;\bx_0)$ for $\bx \in B_{r_0}(\bx_0)$, while $\tilde{R}_\alpha(\bx;\bx_0) = G_\alpha(\bx;\bx_0)$ for $\bx \in \Omega \setminus B_{r_1}(\bx_0)$. 

Substituting $G_\alpha$ in \eqref{Gdecompose} into \eqref{Geq}, we obtain
\[
\mathcal{A}_\alpha (\chi(\bx-\bx_0) u_0(\bx-\bx_0) + \tilde{R}_\alpha) = -1 + \delta(\bx-\bx_0) \,,
\]
or
\begin{equation} \label{Achiu}
	\mA_\alpha \tilde{R}_\alpha = -1 + \delta(\bx-\bx_0) - \mathcal{A}_\alpha(\chi(\bx-\bx_0) u_0(\bx-\bx_0)) \,.
\end{equation}
It now remains to determine $\mathcal{A}_\alpha (\chi(\bx-\bx_0) u_0(\bx-\bx_0))$. To do so, we appeal to the equivalence between applying $\mA_\alpha$ to a function $f(\bx)$ on $\Omega$ and applying $-(-\Delta)^\alpha$ to the (periodic or even, depending on whether $\mAp$ or $\mAn$ is being considered) extension of $f(\bx)$ onto $\bbR^2$. We thus have that 
\BE \label{Achiu2}
	\mathcal{A}_\alpha (\chi(\bx-\bx_0) u_0(\bx-\bx_0)) = -(-\Delta)^\alpha \sum_{\bm \in \bbZ^{2}}  \chi(\bx-T_\bm(\bx_0)) u_0(\bx-T_\bm(\bx_0)) \,, \qquad \bx \in \Omega \,,
\EE
where $T_\bm$ may be either  $T_\bm^{(p)}$ in \eqref{Tm} or $T_\bm^{(n)}$ in \eqref{Y} depending on whether the Green's function being computed is $G_\alpha^{(p)}$ or $G_\alpha^{(n)}$, respectively. Let us consider first the $\bm \neq \bzero$ terms of the sum in \eqref{Achiu2}, given by
\begin{multline} \label{mnonzero}
	\rho_\bm(\bx) \equiv -(-\Delta)^\alpha [ \chi(\bx-T_\bm(\bx_0)) u_0(\bx-T_\bm(\bx_0))] = \\ C_\alpha \int_{\bbR^2} \! \frac{\chi(\by-T_\bm(\bx_0)) u_0(\by-T_\bm(\bx_0)) - \chi(\bx-T_\bm(\bx_0)) u_0(\bx-T_\bm(\bx_0))}{|\by-\bx|^{2+2\alpha}} \, d\by \,, \qquad \bx \in \Omega \,.
\end{multline}
Since $r_1$ in the cut-off function $\chi$ is chosen so that the support of $\chi(\bx-T_\bm(\bx_0))$ in \eqref{mnonzero} lies entirely inside $\Omega_\bm$, we have that $\chi(\bx-T_\bm(\bx_0)) = 0$ for $\bx \in \Omega$. So, \eqref{mnonzero} becomes
\begin{equation} \label{mnonzero2}
	\rho_\bm(\bx) = \\ C_\alpha \int_{|\by| \leq B_{r_1}(T_\bm(\bx_0))} \! \frac{\chi(\by-T_\bm(\bx_0)) u_0(\by-T_\bm(\bx_0))}{|\by-\bx|^{2+2\alpha}} \, d\by \,, \quad \bx \in \Omega \,.
\end{equation}
Since $\bx \in \Omega$ and $\by \in B_{r_1}(T_\bm(\bx_0))$, which has no intersection with $\Omega$, the domain of integration in \eqref{mnonzero2} does not include $\bx$. Thus, the integrand in \eqref{mnonzero2} contains only one singularity, occurring when $\by = T_\bm(\bx_0)$, proportional to $|\by-T_\bm(\bx_0)|^{2\alpha-2}$, which is integrable in $\bbR^2$ when $\alpha \in (0,1)$. The $\bm \neq 0$ terms in \eqref{Achiu2} are thus easily obtained using standard methods of numerical integration.

We now consider the $\bm = \bzero$ term in \eqref{Achiu2},
\begin{multline} \label{mzero}
	-(-\Delta)^\alpha [ \chi(\bx-\bx_0) u_0(\bx-\bx_0)] = \\ C_\alpha \int_{\bbR^2} \! \frac{\chi(\by-\bx_0) u_0(\by-\bx_0) - \chi(\bx-\bx_0) u_0(\bx-\bx_0)}{|\by-\bx|^{2+2\alpha}} \, d\by \,, \qquad \bx \in \Omega \,.
\end{multline}
The integral on the right-hand side of \eqref{mzero} is unbounded at $\bx = \bx_0$ owing to when the variable of integration $\by \to \bx_0$. To extract this singularity, we rewrite the right-hand as
\bes
\begin{multline} \label{mzero2}
C_\alpha \int_{\bbR^2} \! \frac{\chi(\by-\bx_0) u_0(\by-\bx_0) - \chi(\bx-\bx_0) u_0(\bx-\bx_0)}{|\by-\bx|^{2+2\alpha}} \, d\by =  -\chi(\bx-\bx_0) (-\Delta)^\alpha  u_0(\bx-\bx_0) + \rho_0(\bx) \,;
\end{multline}
\BE \label{remainder}
 \rho_\bzero(\bx) \equiv C_\alpha \int_{\bbR^2} \! \frac{u_0(\by-\bx_0)(\chi(\by-\bx_0) - \chi(\bx-\bx_0))}{|\by-\bx|^{2+2\alpha}} \,d\by \,,
\EE
\ees
By \eqref{Gfree}, the first integral on the right-hand side of \eqref{mzero2} is equal to $\chi(\bx-\bx_0)\delta(\bx-\bx_0)$, or simply $\delta(\bx-\bx_0)$, which cancels the $\delta(\bx-\bx_0)$ term in \eqref{Achiu}. 

We now show that $\rho_\bzero(\bx)$ in \eqref{remainder} is finite when $\bx = \bx_0$, where $\bx_0$ is the location of the singularity in $G_\alpha(\bx;\bx_0)$. When $\bx = \bx_0$, the $\chi(\by-\bx_0) - \chi(\bx-\bx_0)$ term in the numerator of \eqref{remainder} vanishes for all $\by \in B_{r_0}(\bx_0)$. The $\mO(|\by-\bx_0|^{-4})$ singularity at $\by = \bx_0$ from the  $\mO(|\by-\bx_0|^{-2+2\alpha})$ and $\mO(|\by-\bx_0|^{-2-2\alpha})$ contributions of $u_0(\by-\bx_0)$ and $|\by-\bx|^{-2-2\alpha}$, respectively, thus need not be considered since the integrand vanishes for all $\by \in B_{r_0}(\bx_0)$. The integrand is finite for all $\by \in \mathbb{R}^2 \setminus B_{r_0}(\bx_0)$, and has decay rate $\sim |\by|^{-4}$ as $|\by| \to \infty$, which is integrable in $\mathbb{R}^2$. Therefore,  $\rho_\bzero(\bx_0)$ is finite.

When $\bx \in B_{r_0}(\bx_0) \setminus \lbrace \bx_0 \rbrace$, the integrand still vanishes for all $\by \in B_{r_0}(\bx_0)$ since then $\chi(\by-\bx_0) = \chi(\bx-\bx_0) = 1$. As in the $\bx = \bx_0$ case, the region of integration is thus $\by \in \mathbb{R}^2\setminus B_{r_0}(\bx_0)$. Both the singularity at $\by = \bx_0$ from $u_0(\by-\bx_0)$ and that at $\by = \bx$ from $|\by-\bx|^{-2-2\alpha}$ appear in $B_{r_0}(\bx_0)$, and so the integrand is always finite in the domain of integration. Moreover, as above, the integrand decays as $\sim |\by|^{-4}$ as $|\by| \to \infty$, which is integrable in $\mathbb{R}^2$. Thus, $\rho_\bzero(\bx)$ is finite when $\bx \in B_{r_0}(\bx_0) \setminus \lbrace \bx_0 \rbrace$.

When $\bx \in \bbR^2\setminus B_{r_1}(\bx_0)$, i.e., outside the support of the cut-off function, we have $\chi(\bx-\bx_0) = 0$. Since $\chi(\by-\bx_0)$ also vanishes when $\by \in \bbR^2\setminus B_{r_1}(\bx_0)$, the domain of integration is simply $\by \in B_{r_1}(\bx_0)$. The $\by = \bx$ singularity from the $|\by-\bx|^{-2-2\alpha}$ term thus lies outside of the domain of integration. The $\mO(|\by-\bx_0|^{-2+2\alpha})$ singularity at $\by = \bx_0$ from $u_(\bx-\bx_0)$ is in the domain of integration, but is integrable when $0<\alpha<1$. The integrand is finite otherwise, and since the domain of integration is also finite, $\rho_\bzero(\bx)$ is finite when $\bx \in \bbR^2\setminus B_{r_1}(\bx_0)$.

When $r_0 \leq |\bx - \bx_0| \leq r_1$, both the $\by = \bx_0$ and $\by = \bx$ singularities lie in the domain of integration. The $\mO(|\by-\bx_0|^{-2+2\alpha})$ singularity at $\by = \bx_0$ is integrable in $\mathbb{R}^2$ when $0<\alpha<1$, while the integral near the $\mO(|\by-\bx|^{-2-2\alpha})$ singularity is well-defined and finite as a principal value since $\chi(\by-\bx_0)$ is at least twice differentiable at $\bx$ (see \cite{valdinoci2009long}). As before, the integrand has decay rate $\sim |\by|^{-4}$ as $|\by| \to \infty$, which is integrable in $\mathbb{R}^2$. Thus, $\rho_\bzero(\bx)$ is finite when $r_0 \leq |\bx - \bx_0| \leq r_1$.

From \eqref{Achiu} - \eqref{remainder}, we thus have a ``smooth'' problem for $\tilde{R}_\alpha$ (i.e., we have replaced the singular $\delta(\bx-\bx_0)$ term with one that is bounded on $\Omega$),
\BE \label{Ralapheq}
	\mA_\alpha \tilde{R}_\alpha(\bx;\bx_0) = -1 - \sum_{\bm\in\bbZ^2} \rho_\bm(\bx)\,; \qquad \bx \in \Omega \,; \qquad \int_\Omega \! \tilde{R}_\alpha(\bx;\bx_0) \,d\by = -\int_\Omega\!\chi(\bx-\bx_0)u_0(\bx-\bx_0) \, d\by \,.
\EE
Note that, by construction, the right-hand side of \eqref{Ralapheq} must be orthogonal to 1. Subject to the integral constraint, \eqref{Ralapheq} yields a unique solution for $\tilde{R}_\alpha$.

In Figs. \eqref{fig:periodicplots} and \eqref{fig:neumannplots} below, we plot $G_\alpha^{(p)}$ and $G_\alpha^{(n)}$ constructed using the above procedure, where we truncate the sum in \eqref{Ralapheq}, \eqref{fraclapsum4}, and \eqref{Ap} at a suitable value of $|\bm|_\infty = m_{\textrm{max}}$ such that changing $m_{\textrm{max}}$ does not appreciably change the result. In Figs. \eqref{fig:rhs_periodic} and \eqref{fig:rhs_neumann}, we plot the right-hand side of \eqref{Ralapheq}. In Figs. \eqref{fig:R_tilde_periodic} and \eqref{fig:R_tilde_neumann}, we plot $\tilde{R}_\alpha^{(p)}(\bx;\bx_0)$ and $\tilde{R}_\alpha^{(n)}(\bx;\bx_0)$ obtained by solving \eqref{Ralapheq} numerically subject to the integral constraint. In Figs. \eqref{fig:G_periodic} and \eqref{fig:G_neumann}, we plot the Green's functions $G_\alpha^{(p)}(\bx;\bx_0)$ and $G_\alpha^{(n)}(\bx;\bx_0)$ satisfying \eqref{Gperiodic} and \eqref{Gneumann}, respectively. Note that the singularity at $\bx_0$ has the exact construction given by \eqref{Gdecompose}.

\begin{figure}
	\centering
	\begin{subfigure}{0.49\textwidth}
		\centering
		\includegraphics[width=1.1\textwidth]{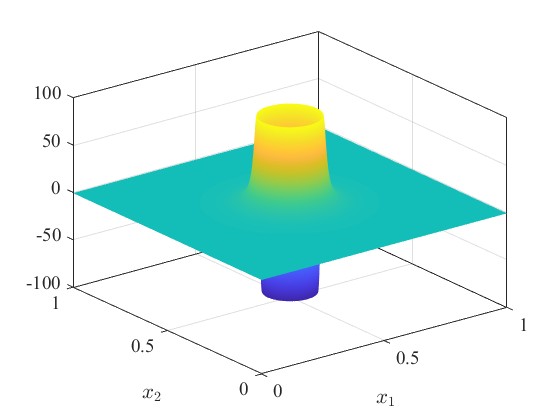}
		\caption{$-1 - \sum_{\bm\in\bbZ^2} \rho_\bm(\bx)$ for periodic BCs}
		\label{fig:rhs_periodic}
	\end{subfigure}
	\begin{subfigure}{0.49\textwidth}
		\centering
		\includegraphics[width=1.1\textwidth]{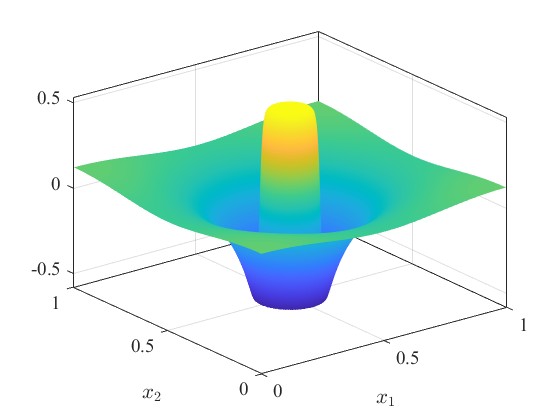}
		\caption{$\tilde{R}_\alpha^{(p)}(\bx;\bx_0)$}
		\label{fig:R_tilde_periodic}
	\end{subfigure}
		\begin{subfigure}{0.49\textwidth}
		\centering
		\includegraphics[width=1.1\textwidth]{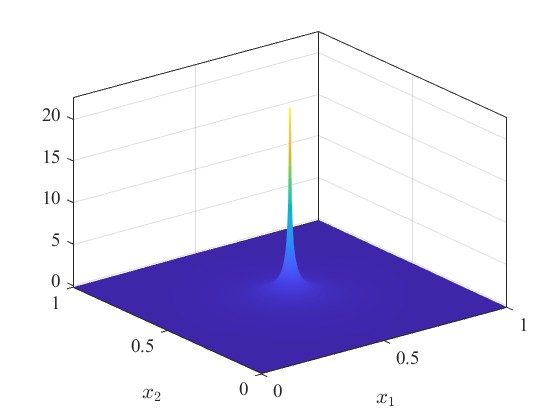}
		\caption{$-G_\alpha^{(p)}(\bx;\bx_0)$}
		\label{fig:G_periodic}
	\end{subfigure}
	\begin{subfigure}{0.49\textwidth}
		\centering
		\includegraphics[width=1.1\textwidth]{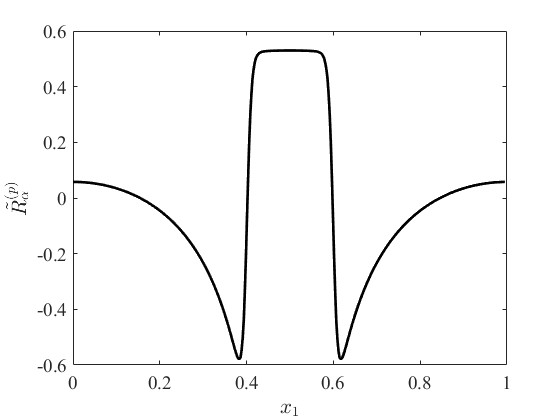}
		\caption{cross section of $\tilde{R}_\alpha^{(p)}(\bx;\bx_0)$}
		\label{fig:R_tilde_periodic_slice}
	\end{subfigure}
	\caption{(a) The right-hand side of \eqref{Ralapheq}, which is orthogonal to 1. (b) The solution of \eqref{Ralapheq}, $\tilde{R}_\alpha^{(p)}(\bx;\bx_0)$. (c) $-G^{(p)}$ as obtained from $\tilde{R}_\alpha^{(p)}$ along with \eqref{Gdecompose}. (d) Cross-section of $\tilde{R}_\alpha^{(p)}(\bx;\bx_0)$ along $x_1$ for $x_2 = 0.5$. Note that the cross-section exhibits periodicity at $x_1 = 0$ and $x_1 = 1$. Here, $\alpha = 0.6$ and $\bx_0 = (0.5, 0.5)$. Note that the gradient of $\tilde{R}_\alpha^{(p)}$ at the location of the singularity, $\bx_0$, is equal to zero.}
	\label{fig:periodicplots}
\end{figure}

 Lastly in Figs. \eqref{fig:R_tilde_periodic_slice} and \eqref{fig:R_tilde_neumann_slice}, we plot cross-sections of $\tilde{R}_\alpha^{(p)}(\bx;\bx_0)$ and $\tilde{R}_\alpha^{(n)}(\bx;\bx_0)$ through $\bx_0$. We recall that $\tilde{R}_\alpha(\bx;\bx_0) = R_\alpha(\bx;\bx_0)$ for $\bx \in B_{r_0}(\bx_0)$, so the correction term of \eqref{ubar}, for example, can be directly obtained from $\tilde{R}_\alpha(\bx;\bx_0)$. We also recall that $\tilde{R}_\alpha(\bx;\bx_0) = G_\alpha(\bx;\bx_0)$ for $\bx \in \Omega \setminus B_{r_1}(\bx_0)$, so the boundary conditions satisfied by $\tilde{R}_\alpha$ are also satisfied by $G_\alpha(\bx;\bx_0)$. In Fig. \ref{fig:R_tilde_periodic_slice}, we observe that $\tilde{R}_\alpha^{(p)}(\bx;\bx_0)$ is periodic at the boundary (in fact, since $\bx_0$ is centered in $\Omega$, it is also Neumann by symmetry); In Fig. \ref{fig:R_tilde_periodic_slice}, $\tilde{R}_\alpha^{(n)}(\bx;\bx_0)$ is Neumann at the boundary despite $\bx_0$ being uncentered in $\Omega$. We note that the quantity $\nabla_{\bx} R_\alpha^{(n)}(\bx;\bx_0)\mid_{\bx = \bx_0}$ is also immediately available from $\tilde{R}_\alpha^{(n)}(\bx;\bx_0)$ -- the first component of $\nabla_{\bx} R_\alpha^{(n)}(\bx;\bx_0)\mid_{\bx = \bx_0}$ is simply the slope at $x_1 = 0.2143$ of Fig. \ref{fig:R_tilde_neumann_slice}. This gradient of the regular part evaluated at the location of the singularity is a quantity that can be used in the context of determining the slow dynamics of localized spot solutions in singularly perturbed reaction-diffusion systems (see, e.g., \cite{ward2002dynamics, kolokolnikov2009spot, ChenWard, gomez2021asymptotic, tzou2019spot, tzou2020analysis}).

\begin{figure}
	\centering
	\begin{subfigure}{0.49\textwidth}
		\centering
		\includegraphics[width=1.1\textwidth]{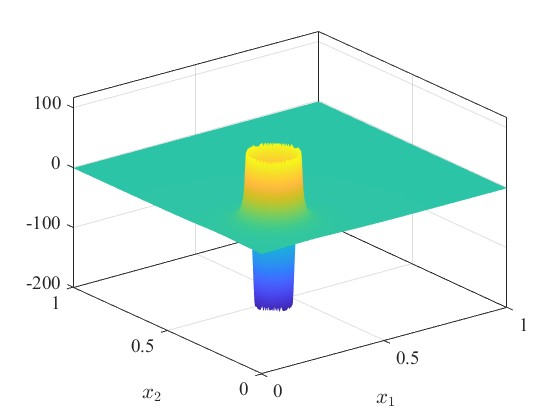}
		\caption{$-1 - \sum_{\bm\in\bbZ^2} \rho_\bm(\bx)$ for Neumann BCs}
		\label{fig:rhs_neumann}
	\end{subfigure}
	\begin{subfigure}{0.49\textwidth}
		\centering
		\includegraphics[width=1.1\textwidth]{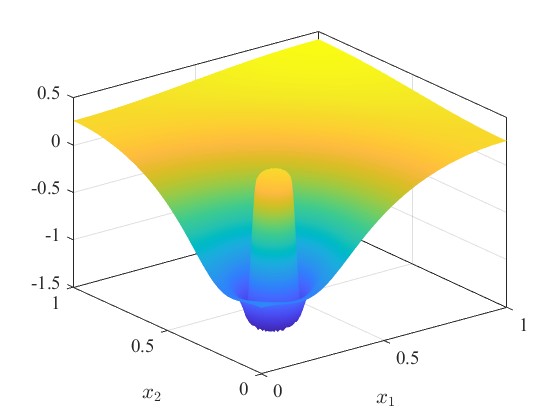}
		\caption{$\tilde{R}_\alpha^{(n)}(\bx;\bx_0)$}
		\label{fig:R_tilde_neumann}
	\end{subfigure}
	\begin{subfigure}{0.49\textwidth}
		\centering
		\includegraphics[width=1.1\textwidth]{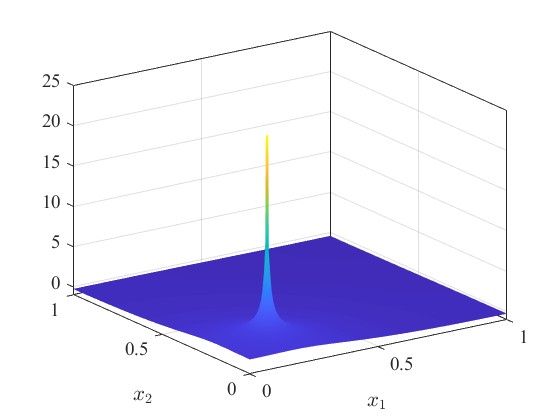}
		\caption{$-G_\alpha^{(n)}(\bx;\bx_0)$}
		\label{fig:G_neumann}
	\end{subfigure}
	\begin{subfigure}{0.49\textwidth}
		\centering
		\includegraphics[width=1.1\textwidth]{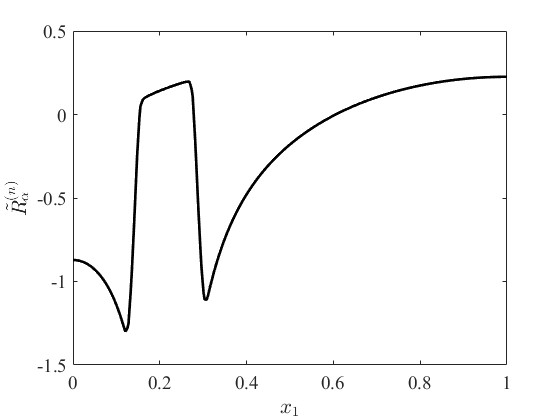}
		\caption{cross section of $\tilde{R}_\alpha^{(n)}(\bx;\bx_0)$}
		\label{fig:R_tilde_neumann_slice}
	\end{subfigure}
	\caption{(a) The right-hand side of \eqref{Ralapheq}, which is orthogonal to 1. (b) The solution of \eqref{Ralapheq}, $\tilde{R}_\alpha^{(n)}(\bx;\bx_0)$. (c) $-G^{(n)}(\bx;\bx_0)$ as obtained from $\tilde{R}_\alpha^{(n)}$ along with \eqref{Gdecompose}. (d) Cross-section of $\tilde{R}_\alpha^{(n)}(\bx;\bx_0)$ along $x_1$ for $x_2 = 0.2143$. Note that the cross-section has zero derivative at $x_1 = 0$ and $x_1 = 1$. Here, $\alpha = 0.6$ and $\bx_0 = (0.2143, 0.2143)$. Note that the gradient of $\tilde{R}_\alpha^{(n)}$ at the location of the singularity, $\bx_0$, is nonzero, owing to boundary effects.}
	\label{fig:neumannplots}
\end{figure}

\section{Application 1: Effects of target configuration and location in the L\'{e}vy flight narrow capture problem} \label{sec:narrow}

In this section, we discuss how the GMFPT for the narrow capture problem is impacted by (1) the placement of a single target within a domain with reflective boundary conditions (see Fig. \ref{fig:schematic_onetrap}), and (2) the relative placement of multiple targets within a domain with periodic boundary conditions (see Fig. \ref{fig:schematic_twotraps}). For (1), we note that a two-term expansion for the GMFPT of a single target within a periodic domain derived in \cite{tzou2024counterexample} and given in \eqref{ubar}. A similar procedure yields the identical expansion for $\bar{u}_\eps^{(n)}$ in the case of Neumann boundary conditions, except with the regular part of the Neumann Green's function, $R_\alpha^{(n)}(\bx_0; \bx_0)$, in the place of its periodic counterpart $R_\alpha^{(p)}(\bx_0;\bx_0)$:
\BE \label{ubarneumann}
	\bar{u}_\eps^{(n)} \sim \eps^{2\alpha-2}\frac{\Gamma(1-\alpha)}{4^\alpha\pi\Gamma(\alpha)} \chi_\alpha - R_\alpha^{(n)}(\bx_0;\bx_0) \,; \qquad \chi_\alpha \equiv \frac{\pi(1-\alpha)}{\sin[(1-\alpha)\pi)]} \,.
\EE

While the leading order terms of $\bar{u}_\eps^{(p)}$ and $\bar{u}_\eps^{(n)}$ are the same, their $\mO(1)$ correction terms behave quite differently. By symmetry, $R_\alpha^{(p)}(\bx_0; \bx_0)$ is independent of $\bx_0$, and so the GMFPT $\bar{u}_\eps^{(p)}$ remains unchanged no matter where in $\Omega$ the target is centered. This not true, however, in the Neumann case where boundary effects are important. This is shown in Fig. \ref{fig:mfpt_vs_trap_location_asympt_vs_numerical_neumann}, which shows how $\bar{u}_\eps^{(n)}$ varies as a function of the target's position in $\Omega$, $(s,s)$ for $0 < s \leq 0.5$. The red curve is obtained by numerically solving \eqref{narrowescapeequation} via finite differences, while the blue is obtained from \eqref{ubarneumann} with $R_\alpha^{(n)}(\bx_0;\bx_0)$ computed from \S \ref{sec:Greens}. All variation with respect to the variable $s$ is due to the variation of $R_\alpha^{(n)}(\bx_0; \bx_0)$ as $\bx_0$ is varied in $\Omega$. By symmetry, $\bar{u}_\eps^{(n)}$ reaches its minimum when the target is at the center of $\Omega$ (i.e., when $s = 1/2$). When the target is near $\partial \Omega$ (e.g., when $s$ is small), the target is partially ``shielded'' by the boundary, contributing to a higher GMFPT. We discuss shielding effects in more detail in \S \ref{sec:split}.

Before deriving the expansion for $\bar{u}_\eps^{(p)}$ for multiple targets, we highlight the utility of our hybrid asymptotic-numerical approach in contrast to the straightforward numerical solve of \eqref{narrowescapeequation} or its periodic counterpart. We begin by noting that the nonlocal nature of $\mAp$ and $\mAn$ means that their discretization results in a full matrix, $A_\alpha$. This contrasts the sparse matrix that results from discretization of the regular Laplacian operator. Assuming a uniform grid spacing of $h$, $A_\alpha$ must then have $\mO(1/h^4)$ entries. Each of these entries involves computing the distance from one grid point to every other grid point (as well as their corresponding image points). When the target size $\eps$ is small, the computational time it takes to populate the matrix $A_\alpha$, along with the memory needed to store it, can become prohibitive. Because $h$ must scale in proportion to $\eps$, the computational time and memory requirements scale must then scale as $\mO(1/\eps^{4})$ (symmetry in the periodic case can be used to circumvent the $\mO(1/\eps^{4})$ for the computational time).

In our hybrid asymptotic-numerical approach, the small $\eps$ scale is handled asymptotically, leaving the numerical problem for $G_\alpha$ free of $\eps$. Furthermore, by removing the singular part of $G_\alpha$ analytically, we are left with a numerical problem without features requiring a small grid spacing to resolve. Thus, this numerical problem can employ a grid spacing $h$ that is larger than that required for the $\eps$-dependent problem - a doubling of the grid spacing can result in a 16-times reduction in the computational time and memory requirements. We should note, however, that because the support of the cutoff function must lie entirely inside $\Omega$ (otherwise $\tilde{R}_\alpha$ would not satisfy the same boundary conditions as $G_\alpha$), the numerical method of \S \ref{sec:Greens} may require a smaller grid spacing when $\bx_0$ is very near the boundary due to the large gradient of $\chi$ in \eqref{cutoff}.

We now derive a two-term expansion for $\bar{u}_\eps$ using a matched asymptotic analysis. The procedure is the same for $\bar{u}_\eps^{(p)}$ and $\bar{u}_\eps^{(n)}$, and hence we will dispense with the superscripts as in \S \ref{sec:Greens}. The method follows that of, e.g., \cite{gomez2024first, ward1993summing, cheviakov2010asymptotic, pillay2010asymptotic, kolokolnikov2005optimizing, tzou2024counterexample}, and will thus be brief. Let $u_\eps(\bx)$ denote the mean first passage time to one of a set of $N$ targets of radius $\eps_i \equiv \kappa_i\eps$ centered at $\bx_i\in \Omega$. Here, $0<\eps \ll 1$, $\kappa_i > 0$ and $\kappa_i \sim \mO(1)$ with respect to $\eps$ for $i = 1, \ldots, N$. We assume that the target centers are well-separated, i.e., $|\bx_i - \bx_j| \sim \mO(1)$ with respect to $\eps$ for $i \neq j$. Then $u_\eps(\bx)$ satisfies
\bes \label{narrowescapeequationmult}
\BE \label{narrowescapeequationmulteq}
\mA u_\eps(\bx) = -1 \,, \quad  \bx \in \Omega \setminus \bigcup_{i = 1}^N  B_{\eps_i}(\bx_i) \,; \qquad  u_\eps = 0 \enspace \mbox{for} \enspace \bx \in \bigcup_{i = 1}^N  B_{\eps_i}(\bx_i) \,; \quad \bx_i \in \Omega \enspace \mbox{for} \enspace i = 1, \ldots, N \,.
\EE
\BE \label{narrowescapeequationmultec}
	u_\eps = 0 \enspace \mbox{for} \enspace \bx \in \bigcup_{i = 1}^N  B_{\eps_i}(\bx_i) \,; \quad \bx_i \in \Omega \enspace \mbox{for} \enspace i = 1, \ldots, N \,; \qquad \eps_i \equiv \kappa_i \eps \,, \quad i = 1,\ldots,N \,.
\EE
\ees
The GMFPT, $\bar{u}_\eps$, is the average of $u_\eps$ over $\Omega$,
\BE \label{ubarformula}
	\bar{u}_\eps = \frac{1}{\Omega}\int_\Omega \! u_\eps(\bx) \,, d\bx \,.
\EE
We begin in the $\mO(\eps_i)$ region centered at $\bx_i$, and let
\bes \label{innarvar}
\BE \label{changeofvar}
	\bx = \bx_i + \eps_i\bz \,; \qquad U_i(\bz) = u_{\eps}(\bx_i + \eps_i \bz) \,.
\EE
We expand the inner variable $U_i$ as
\BE \label{innerexpansion}
	U_i \sim \eps_i^{2\alpha-2} U_{i0} + \eps_i^{2\alpha} U_{i1} \,.
\EE
\ees
As shown in \cite{tzou2024counterexample}, $\mA_\alpha f(\eps_i^{-1}\bx) \sim -\eps_i^{-2\alpha}(-\Delta_\mathbf{z})^\alpha f(\mathbf{z})$ as $\eps_i \to 0^+$, where $-(-\Delta_\mathbf{z})^\alpha$ is the fractional Laplacian with respect to the variable $\mathbf{z}$.

Substituting \eqref{changeofvar} and using \eqref{innerexpansion}, we obtain the leading order inner equation near $\bx_i$,
\bes
\BE \label{Ui0eq}
	-(-\Delta_\bz)^{\alpha} U_{i0} = 0 \,, \quad \bz \in \bbR^2 \setminus B_{1}(\bzero) \,, \qquad U_{i0} = 0 \,, \quad \bz \in B_1(\bzero) \,,
\EE
where the target that is the disk of radius $\eps_i$ centered about $\bx_i$ is now the unit disk centered at the origin due to the change of variables \eqref{changeofvar}. As $|\bz| \to \infty$, we assume that $U_{i0}$ is radially symmetric and takes the form
\BE \label{Ui0far}
	U_{i0} \sim S_i \left( -\frac{1}{|\bz|^{2-2\alpha}} + \chi_\alpha   \right) \enspace \mbox{as} \enspace |\bz| \to \infty \,.
\EE
\ees
In \eqref{Ui0far}, $\chi_\alpha$ is an $\mO(1)$ constant given in \eqref{ubar} that depends on $\alpha$ as well as the geometry of the rescaled target. Since we have assumed that all targets are in the form of a circular disk, for which an explicit formula exists for $\chi_\alpha$ (see \cite{tzou2024counterexample} and \cite{kahane1981solution}), $\chi_\alpha$ is the same for each inner region. For general target geometries, a numerical solution of a certain integral equation may be required to obtain $\chi_\alpha$.

In the outer region, in the limit $\eps \to 0^+$, the exterior condition near $\bx_i$ of \eqref{narrowescapeequationmult} is replaced by a local behavior, obtained from the far-field behavior of $U_{i0}$ in \eqref{Ui0far} with \eqref{innarvar}, that specifies both the singular structure of $u_\eps$ as $\bx \to \bx_i$, as well as the regular part, yielding
\bes \label{narrowescapeequationmult2}
\BE \label{narrowescapeequationmult2eq}
\mA u_\eps(\bx) = -1 \,, \quad  \bx \in \Omega \setminus \bigcup_{i = 1}^N  \lbrace \bx_i \rbrace \,;
\EE
\BE \label{narrowescapeequationmult2loc}
u_\eps \sim \eps_i^{2\alpha-2} S_i\left(-\frac{\eps_i^{2-2\alpha}}{|\bx-\bx_i|^{2-2\alpha}} + \chi_\alpha \right)  \enspace \mbox{as} \enspace \bx \to \bx_i \,, \quad i = 1, \ldots, N \,.
\EE
\ees
Comparing \eqref{narrowescapeequationmult2} to \eqref{G}, the leading order solution to $u_\eps$ may be written as
\BE \label{ueps}
	u_\eps \sim \frac{1}{c_\alpha}\sum_{i = 1}^N S_i G_\alpha(\bx;\bx_i) + \bar{u}_\eps \,.
\EE
By \eqref{ubarformula} along with the zero-integral condition for $G_\alpha$ in \eqref{Geq}, $\bar{u}_\eps$ in \eqref{ueps} is the mean of $u_\eps$ over $\Omega$; i.e., $\bar{u}_\eps$ is the GMFPT of the first passage process described by \eqref{narrowescapeequationmult}. 

We now formulate a system of $N+1$ equations for $\bar{u}_\eps$ and $S_i$, $i = 1, \ldots, N$. We note first that 
\BE \label{fraclapueps}
	\mA_\alpha u_\eps = \frac{1}{c_\alpha}\sum_{i = 1}^N S_i[-1 + \delta(\bx - \bx_i)] \,.
\EE
Comparing \eqref{fraclapueps} to \eqref{narrowescapeequationmult2} for $\bx \in \Omega \setminus \bigcup_{i = 1}^N \lbrace \bx_i\rbrace$, we require the consistency condition for $S_i$,
\BE \label{sumS}
	\sum_{i = 1}^N S_i = c_\alpha \,.
\EE
To obtain the other $N$ equations, we match the local behavior near $\bx_i$ of $u_\eps$ given by \eqref{ueps} to that required by \eqref{narrowescapeequationmult2loc}. The former is given by
\BE \label{locmatch}
	u_\eps \sim \frac{S_i}{c_\alpha}\left(  -\frac{c_\alpha}{|\bx-\bx_i|^{2-2\alpha}} + R_\alpha(\bx_i;\bx_i) \right) + \frac{1}{c_\alpha}\sum_{j\neq i}^N S_j G_\alpha(\bx_i; \bx_j) + \bar{u}_\eps \,, \qquad i = 1,\ldots, N \,.
\EE
Matching the local behavior of \eqref{locmatch} to that in \eqref{narrowescapeequationmult2loc}, we obtain the $N$ equations
\BE \label{locmatch2}
	S_iR_\alpha(\bx_i;\bx_i) + \sum_{j\neq i}^N S_j G_\alpha(\bx_i; \bx_j) + c_\alpha\bar{u}_\eps = c_\alpha\eps^{2\alpha-2} \kappa_i^{2\alpha-2}S_i\chi_\alpha \,, \qquad i = 1,\ldots, N \,.
\EE
In matrix form, the $N+1$ equations given by \eqref{sumS} together with \eqref{locmatch2} are
\bes\label{Nplus1}
\BE \label{Nplus1eq}
 \mathcal{G}_\alpha \bs + c_\alpha \bar{u}_\eps \be = c_\alpha\chi_\alpha \eps^{2\alpha-2} \mathcal{K} \bs \,; \qquad \be^T\bs = c_\alpha \,,
\EE
where the entries of the $N\times N$ matrices $ \mathcal{G}_\alpha$ and $\mathcal{K}$, along with the $N\times 1$ vectors $\bS$ and $\be$ are given by
\BE \label{Nplus1def}
	 \mathcal{G}_\alpha^{(ij)} = \begin{cases} 
	 	R_\alpha(\bx_i;\bx_i) & i = j \,, \\
	 	G_\alpha(\bx_i;\bx_j) & i \neq j \,,
	 \end{cases};
	 \quad \mathcal{K}^{(ij)} = \begin{cases} 
	 	\kappa_i^{2\alpha-2} & i = j \,, \\
	 	0 & i \neq j \,, 
	 \end{cases}; \quad
	 \be = \begin{pmatrix}
	 	1 \\ \vdots \\ 1
	 \end{pmatrix}; \quad
	 \bs = \begin{pmatrix}
	 	S_1 \\ \vdots \\ S_N
	 \end{pmatrix} \,.
\EE
\ees

Finally, it remains to solve for $\bar{u}_\eps$ in \eqref{Nplus1eq}. First, we invert for $\bs$ to find
\BE \label{ssol}
	\bs = c_\alpha \bar{u}_\eps\left\lbrack c_\alpha\chi_\alpha\eps^{2\alpha-2}\mathcal{K} - \mathcal{G}_\alpha \right\rbrack^{-1}\be \,.
\EE
Note that $c_\alpha\chi_\alpha\eps^{2\alpha-2}\mathcal{K} - \mathcal{G}_\alpha$ must be a strictly diagonally dominant matrix in the limit $\eps \to 0^+$, since $\mathcal{K}$ is diagonal. Hence, it must be invertible. We then take the inner product of both sides of \eqref{ssol} with $\be$, and use the second equation of \eqref{Nplus1eq} to obtain for $\bar{u}_\eps$,
\BE \label{ubarmult}
	\bar{u}_\eps \sim \frac{1}{\be^T \left\lbrack c_\alpha\chi_\alpha\eps^{2\alpha-2}\mathcal{K} - \mathcal{G}_\alpha \right\rbrack^{-1}\be} \,.
\EE
In the case $N=1$ of the single target of radius $\eps$ at $\bx_0 \in\Omega$, we have $\mathcal{G}_\alpha = R_\alpha(\bx_0;\bx_0)$, and \eqref{ubarmult} reduces to \eqref{ubar}. Expanding \eqref{ubarmult} to two orders, we obtain
\BE \label{ubarmult2}
	\bar{u}_\eps \sim \frac{c_\alpha\chi_\alpha\eps^{2\alpha-2}}{\be^T\mathcal{K}^{-1}\be} - \frac{1}{[\be^T \mathcal{K}^{-1} \be]^2}\be^T \mathcal{K}^{-1}\mathcal{G}_\alpha \mathcal{K}^{-1} \be \,.
\EE
In the case where all targets are disks of radius $\eps$, we have $\kappa_i = 1$ for $i = 1,\ldots, N$, and $\mathcal{K} = I$, where $I$ is the $N\times N$ identity matrix.  In this case, \eqref{ubarmult2} simplifies to
\BE \label{ubarmult3}
	\bar{u}_\eps \sim \frac{c_\alpha\chi_\alpha\eps^{2\alpha-2}}{N} - \frac{1}{N^2}\be^T \mathcal{G}_\alpha \be \,.
\EE
Compared to \eqref{ubar}, the leading order term is simply scaled by $N$, the number of targets, while the $\mO(1)$ correction term has been augmented to include the off-diagonal terms of the Green's interaction matrix $\mathcal{G}_\alpha$, accounting for interactions between all possible pairs of targets. It is this correction term that captures the change of $\bar{u}_\eps^{(p)}$ in Fig. \ref{fig:mfpt_vs_second_trap_location_asympt_vs_numerical_periodic} as the parameter $s$ is varied. For $s$ near $(0.25, 0.25)$, the two targets are relatively close together (see Fig. \ref{fig:schematic_twotraps}), inducing a shielding effect on each other. As $s$ is increased and the two targets become more evenly spaced, the shielding effect is lessened, and the GMFPT decreases. By symmetry, the GMFPT is minimized when $s = 0.75$.

\section{Application 2: Splitting probabilities in the Brownian and L\'{e}vy flight narrow capture problems }\label{sec:split}

In this section, we derive a two-term expansion for the splitting probability to reach one particular target inside the domain \textit{before} hitting any of the other targets. See Fig. \ref{fig:schematic_splitting} for a schematic representation. Let us denote by $v_\eps(\bx)$ the probability of hitting a target of radius $\eps_0 = \kappa_0\eps$ centered at $\bx_0$ before hitting any of the $N$ other targets of radius $\eps_i = \kappa_i\eps$ centered at $\bx_i \in \Omega$. As before, we have $0 < \eps \ll 1$ and $\kappa_i \sim \mO(1)$ with respect to $\eps$ for $i = 0,\ldots, N$. We assume that the target centers are well-separated, i.e., $|\bx_i - \bx_j| \sim \mO(1)$ with respect to $\eps$ for $i \neq j$, and also that $\min_{\bx \in \partial\Omega}|\bx_i-\bx| \sim \mO(1)$ with respect to $\eps$ for $i = 0, \ldots, N$. Following a similar derivation to \cite{tzou2024counterexample} of the elliptic problem for the narrow escape time, and in analogy with the Brownian splitting probability \cite{paquin2023narrow, redner2001guide, bressloff2022narrow, delgado2015conditional, kurella2015asymptotic}, $v_\eps(\bx)$ satisfies
\bes \label{splitting}
\BE \label{splittingeq}
	\mA v_\eps = 0 \,, \quad \bx \in \Omega \setminus \bigcup_{i = 0}^N B_{\eps_i}(\bx_i) \,,
\EE
\BE \label{splittingec}
	v_\eps = 1 \enspace \mbox{for} \enspace \bx \in B_{\eps_0}(\bx_0) \,, \qquad v_\eps = 0 \enspace \mbox{for} \enspace \bx \in \bigcup_{i = 1}^N B_{\eps_i}(\bx_i) \,; \qquad \eps_i \equiv \kappa_i \eps \,, \quad i = 0,\ldots,N \,.
\EE
\ees
Note that $v_\eps = 1$ when the starting location is inside the desired target at $\bx_0$, while $v_\eps = 0$ when the starting location is inside an obstacle target at $\bx_i$ for $i = 1,\ldots, N$. We define the average splitting probability, $\bar{v}_\eps$, as the average of $v_\eps$ in \eqref{splitting} over $\Omega$,
\BE \label{vbarformula}
	\bar{v}_\eps = \frac{1}{\Omega} \int_\Omega \! v_\eps(\bx) \, d\bx \,.
\EE

We now derive a two term expansion for $\bar{v}_\eps$. The procedure is similar to that of \S \ref{sec:narrow} for the GMPFT, so what follows will be brief. In the $\mO(\eps)$ region centered at $\bx_i$, $i = 0, 1, \ldots, N$, we let
\BE \label{innearvarssplit}
	\bx = \bx_i + \eps_i \bz \,; \qquad V_i(\bz) = v_\eps(\bx_i + \eps_i \bz) \,.
\EE
Letting $V_i \sim \eps^{2\alpha-2} V_{i0}$, we have for the $i$-th inner region
\bes \label{Vi}
\BE \label{Vi0eq}
-(-\Delta_\bz)^{\alpha} V_{i0} = 0 \,, \quad \bz \in \bbR^2 \setminus B_{1}(\bzero) \,, \quad V_{00} = 1 \enspace \mbox{and} \enspace V_{i0} = 0 \,, \quad \bz \in B_1(\bzero) \,, \quad i = 1,\ldots, N \,,
\EE
supplemented by the far-field behavior
\BE \label{Vifar}
	V_{00} \sim S_0 \left( -\frac{1}{|\bz|^{2-2\alpha}} + \chi_\alpha   \right) + 1 \enspace \mbox{and} \enspace	V_{i0} \sim S_i \left( -\frac{1}{|\bz|^{2-2\alpha}} + \chi_\alpha   \right) \enspace \mbox{as} \enspace |\bz| \to \infty \,, \quad i = 1, \ldots, N \,.
\EE
\ees
In the limit $\eps \to 0^+$ in \eqref{splitting}, we obtain
\bes \label{splittinglim}
\BE \label{splittinglimeq}
	\mA v_\eps = 0 \,, \quad \bx \in \Omega \setminus \bigcup_{i = 0}^N \lbrace \bx_i \rbrace \,,
\EE
along with the required singularity conditions near $\bx_i$, $i = 0, \ldots, N$, from \eqref{Vifar}
\BE \label{splittinglimloci}
v_\eps  \sim \eps_i^{2\alpha-2} S_i \left( -\frac{\eps_i^{2-2\alpha}}{|\bx-\bx_i|^{2-2\alpha}} + \chi_\alpha   \right) + \delta_{0i} \enspace \mbox{as} \enspace \bx \to \bx_i \,, \quad 0 = 1, \ldots, N \,,
\EE
\ees
where $\delta_{0i}$ is the Kronecker delta function.

Comparing \eqref{splittinglim} to \eqref{G}, we may write
\BE \label{veps}
	v_\eps(\bx) \sim \frac{1}{c_\alpha}\sum_{i = 0}^N S_iG_\alpha(\bx;\bx_i) + \bar{v}_\eps \,,
\EE
where the weights $S_i$ must satisfy the consistency condition
\bes \label{splitequations}
\BE \label{sumSsplit}
	\sum_{i = 0}^N S_i = 0 \,,
\EE
owing to the fact that the right-hand side of \eqref{splittinglimeq} is homogeneous. To find the other $N+1$ linear equations for $S_i$, $i = 0,\ldots, N$ and $\bar{v}_\eps$, we match the singularity condition of $v_\eps$ in \eqref{veps} to that required by \eqref{splittinglimloci}, obtaining
\BE \label{splitmatch}
	\frac{S_i}{c_\alpha}R_\alpha(\bx_i;\bx_i) + \frac{1}{c_\alpha} \sum\limits_{\substack{j = 0 \\ j \neq i}}^N S_j G_\alpha(\bx_i;\bx_j) + \bar{v}_\eps = \eps^{2\alpha-2}\kappa_i^{2\alpha-2} \chi_\alpha S_i + \delta_{0i} \,, \quad i = 1\ldots,N \,.
\EE
\ees
In matrix vector form, \eqref{splitequations} becomes
\BE \label{Nplus2spliteq}
\mathcal{G}_\alpha \bs + c_\alpha \bar{v}_\eps \be = c_\alpha\chi_\alpha \eps^{2\alpha-2} \mathcal{K} \bs + c_\alpha\be_1 \,; \qquad \be^T\bs =0 \,,
\EE
where $\mathcal{G}_\alpha$, $\bs$, and $\mathcal{K}$ are given in \eqref{Nplus1def}, except they have dimension $N+1$ in \eqref{Nplus2spliteq} with a starting index of 0 instead of 1, and $\be_1$ is the $(N+1)$-vector given by $\be_1 \equiv (1, 0, \ldots, 0)^T$.

From \eqref{Nplus2spliteq}, $\bs$ is given by
\BE \label{s}
	\bs = c_\alpha\left[ c_\alpha \chi_\alpha \eps^{2\alpha-2} \mathcal{K} - \mathcal{G}_\alpha \right]^{-1}(\bar{v}_\eps \be - \be_1) \,.
\EE
Applying the zero-sum condition on $\bs$ from \eqref{Nplus2spliteq}, we obtain that $\bar{v}_\eps$ is given by
\BE \label{vbar}
	\bar{v}_\eps = \frac{\be^T\left[ c_\alpha \chi_\alpha \eps^{2\alpha-2} \mathcal{K} - \mathcal{G}_\alpha \right]^{-1}\be_1}{\be^T\left[ c_\alpha \chi_\alpha \eps^{2\alpha-2} \mathcal{K} - \mathcal{G}_\alpha \right]^{-1}\be} \,.
\EE
The formula for $\bar{v}_\eps$ in \eqref{vbar} effectively sums all powers of $\eps^{2-2\alpha}$. For a more informative formula for $\bar{v}_\eps$, we can expand \eqref{vbar} for small $\eps$, which yields,
\BE \label{vbar2}
	\bar{v}_\eps \sim \frac{\be^T\mK^{-1}\be_1}{\be^T\mK^{-1}\be} + \frac{\eps^{2-2\alpha}}{c_\alpha\chi_\alpha \left(\be^T \mK^{-1}\be\right)^2}\be^T \mK^{-1} \mG_\alpha \mK^{-1}\left[\left(\be^T\mK^{-1}\be\right) \be_1 - \left(\be^T\mK^{-1}\be_1\right) \be \right] \,.
\EE
And if all targets have radius $\eps$ (i.e., $\mK = I$, where $I$ here is the $(N+1)\times (N+1)$ identity matrix), \eqref{vbar2} simplifies to
\BE \label{vbaridentical}
	\bar{v}_\eps \sim \frac{1}{N+1} + \frac{\eps^{2-2\alpha}}{c_\alpha\chi_\alpha(N+1)^2} \be^T\mG_\alpha \left[ (N+1) \be_1 - \be \right]  \,.	
\EE
The leading order term in \eqref{vbaridentical} is a function only of the number of ``desired'' targets (in this case, we have assumed one) and the number of obstacle targets ($N$). The $\mO(\eps^{2-2\alpha})$ correction term accounts for how the targets are located within the domain and their positions relative to one another.

In Fig. \ref{fig:schematic_splitting}, the desired target (center, heavy line) is surrounded, or ``shielded'', by five obstacle targets. All targets have radius $0 < \eps \ll 1$. In Fig. \ref{fig:split_prob_vs_alpha}, we plot the average splitting probability versus the L\'{e}vy flight index $\alpha$. The red curve is obtained from a finite difference solution of \eqref{splitting}, while the $\alpha<1$ portion of the blue curve are given by  $\bar{v}_\eps$ in \eqref{vbar}, with the Green's functions computed by the methods of \S \ref{sec:Greens}. The $\alpha = 1$ asymptotic estimate comes from the Brownian analog of \eqref{vbar} (not given), the derivation of which is similar to that leading to \eqref{vbar}.

With $N = 5$, in the absence of the effect of relative positions of desired and obstacle targets, the leading order of $\bar{v}_\eps$ given by \eqref{vbaridentical} simply yields $1/6$. In Fig. \ref{fig:split_prob_vs_alpha}, we see that $\bar{v}_\eps \sim 0.165$ when $\alpha = 0.2$, meaning that relative positions play little role when $\alpha$ is small. In this instance, the frequency of long jumps is relatively high; in combination with the fact that the derivation leading to \eqref{splitting} allows jumps over obstacles (see also \cite{gomez2024first, tzou2024counterexample}), the near-absence of ``shielding effects'' (see, e.g., \cite{kurella2015asymptotic}) in the small-$\alpha$ case is not surprising. This is also seen clearly in \eqref{vbaridentical}, where the $\mO(\eps^{2-2\alpha})$ scaling of the correction term means that it becomes asymptotically smaller as $\alpha$ becomes smaller. By contrast, Fig. \ref{fig:split_prob_vs_alpha} shows that $\bar{v}_\eps$ decreases to  $\sim 0.079$ in the $\alpha = 1$ (Brownian) limit.

\begin{figure}
	\centering
		\begin{subfigure}{0.49\textwidth}
		\centering
		\includegraphics[width=1.1\textwidth]{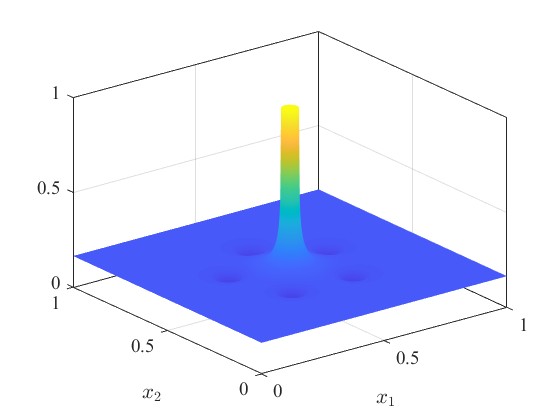}
		\caption{$v_\eps(\bx)$ for $\alpha = 0.2$}
		\label{fig:split_prob_surf_plot_RUN_1_alpha_0p2}
	\end{subfigure}
	\begin{subfigure}{0.49\textwidth}
		\centering
		\includegraphics[width=1.1\textwidth]{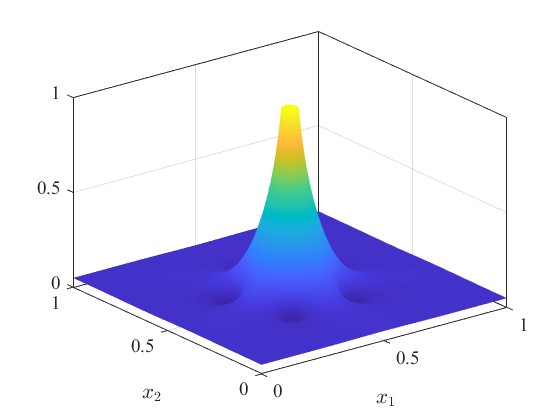}
		\caption{$v_\eps(\bx)$ for $\alpha = 1$}
		\label{fig:split_prob_surf_plot_laplace}
	\end{subfigure}
	\caption{Numerical solution of \eqref{splitting} for the splitting probability $v_\eps(\bx)$ corresponding to the configuration in Fig. \ref{fig:schematic_splitting}. Here, the common radius of all targets is $\eps = 0.03$ while $\alpha = 0.2$ in (a) and $\alpha = 1$ in (b). The ``desired' target is centered at $(0.5, 0.5)$, while the ``obstacle'' targets are located at $(0.5,0.5) + 7\eps(\cos\theta, \sin\theta)$ for $\theta = \pi/2 + 2n\pi/5$ for $n = 0, 1, \ldots, 4$. The decay from $v_\eps = 1$ in (a) is faster than in (b) owing to the power-law decay of the Green's function for $\alpha < 1$ versus the logarithmic behavior for $\alpha = 1$. Outside the ring of obstacle targets, however, $v_\eps$ is closer to 0 in (b) due to the likelihood of a particle hitting an obstacle target first when starting outside of the ring of obstacles. In effect, the obstacle targets ``shield'' the desired target from particles that start outside of the ring. In contrast, particles in (a) undergo long jumps (which can pass over the obstacle targets) with relatively high frequency, mitigating the shielding effect.}
	\label{fig:splitting}
\end{figure}

In Fig. \ref{splitting}, we contrast $v_\eps(\bx)$ for $\alpha = 0.2$ (Fig. \ref{fig:split_prob_surf_plot_RUN_1_alpha_0p2}) and the Brownian limit $\alpha = 1$ (Fig. \ref{fig:split_prob_surf_plot_laplace}) in the case where the desired target is centered at $(0.5, 0.5)$ while the obstacle targets are located at $(0.5,0.5) + 7\eps(\cos\theta, \sin\theta)$ for $\theta = \pi/2 + 2n\pi/5$ for $n = 0, 1, \ldots, 4$ with $\eps = 0.03$. Notice that $v_\eps(\bx) = 1$ when $|\bx-(0.5,0.5)| < \eps$, and $v_\eps(\bx) = 0$ at the locations of the obstacle targets.  

From \eqref{veps} and the Green's function of the fractional Laplacian \eqref{Gloc}, $v_\eps$ in Fig. \ref{fig:split_prob_surf_plot_RUN_1_alpha_0p2} ($\alpha = 0.2$) follows a power-law decay away from $(0.5,0.5)$, which is steeper than that seen in Fig. \ref{fig:split_prob_surf_plot_laplace} ($\alpha = 1$), which exhibits logarithmic decay due to the corresponding Green's function of the Laplacian. This is due to the likelihood of particles that start near $(0.5, 0.5)$ taking a long jump and effectively having to restart the random walk from a far-away location. However, outside of the ring of obstacle targets, $v_\eps$ is closer to 0 in the $\alpha = 1$ case due to the likelihood of particles starting there to hit one of the obstacles prior to reaching the desired target. This shielding effect is far less impactful for L\'{e}vy flight particles, which are able to take long jumps with relatively high frequency and jump over the obstacle targets.

\section{Discussion} \label{sec:disc}

On the unit square $\Omega$, we have derived an expression for an infinitesimal generator $\mAn$ describing a L\'{e}vy flight process with index $\alpha \in (0,1)$ where we assume specular reflection of particles on $\partial\Omega$. We gave a heuristic argument showing this operator is equivalent to the spectral representation of the fractional Laplacian on $\Omega$ with restriction to the set of Laplacian eigenfunctions with zero normal derivative on $\partial\Omega$. For $\mAn$, as well as its periodic analog $\mAp$, we proposed a method for accurately computing its source-neutral Green's function. In particular, the method accurately determines the value of the regular part of the Green's function at the location of the singularity.

We applied this method within an asymptotic framework to determine narrow capture times for two targets in a periodic square and for one target in an insulating square. For the former, our method accurately determined the effect on the average search time of target configuration -- i.e., how the two targets are placed relative to one another within the search domain. For the latter, our method accurately determined how the average search time is affected by where a single target is placed relative to an reflective boundary. Note that both of these effects arise as $\mO(1)$ correction terms that are not captured by the leading order theory and which require accurate computation of Green's functions to obtain.

Using a similar approach, we analyzed splitting probabilities in the narrow capture problem. In particular, we highlighted the difference between how effective L\'{e}vy flights and Brownian motion are at navigating obstacles to reach a desired target. Our example consisted of a desired target surrounded by five obstacles, for which we computed the average probability of a particle undergoing a L\'{e}vy flight of index $\alpha$ to reach the target before hitting any of the obstacles. Our asymptotic analysis along with our method for computing Green's functions showed that this probability is a decreasing function of $\alpha$. In fact, for sufficiently small $\alpha$, the probability was close to the uniform value of $1/6$ indicating almost no shielding effect. On the other hand, the probability was decreased by approximately a factor of $2$ when $\alpha$ approached its Brownian limit of $1$.

This approach of using matched asymptotic methods to recast an $\eps$-dependent problem into a canonical Green's function problem, which can then be accurately solved numerically, has the benefit of avoiding having to solve a numerical problem where the grid spacing scales poorly with $\eps$. This is especially relevant in the two-dimensional nonlocal problems considered in \S \ref{sec:narrow} and \ref{sec:split}, where discretization of the operator results in a full matrix that may become very memory-intensive to store and invert as the number of grid points is increased. Our methods can also be applied to the three-dimensional analogs, where this issue becomes even more accute.

Another type of problem where this approach may be beneficial is that of determining the stability and slow dynamics of localized spot solutions in singularly perturbed reaction-diffusion systems exhibiting L\'{e}vy flights in bounded domains. The use of asymptotic methods and Green's functions to extract detailed results and stability thresholds in the Brownian case was pioneered in, e.g., \cite{ward2002dynamics, ward2002existence, kolokolnikov2009spot}. Since then, efforts have been made to numerically compute Green's functions on general manifolds in order to extend the theory beyond domains in which Green's functions are explicitly known \cite{tzou2019spot, tzou2020analysis}. More recently, \cite{gomez2023multi,gomez2024spike} have extended this theory to L\'{e}vy flights in one spatial dimension, where the Green's function was computed in terms of a rapidly converging series of eigenfunctions. Our methods of \S \ref{sec:Greens} can facilitate further extending the theory to two spatial dimensions, including those on which explicitly known eigenfunctions are not available.

While the square domain we considered was particularly convenient for our method of images-type approach to computing the Neumann Green's function of $\mAn$, we expect that our approach can be easily adapted to computing the periodic source-neutral Green's functions of $\mAp$ on general periodic Bravais lattices. This then makes it possible to compute the principal eigenvalue of the fractional Laplacian on such lattices. The was done for the Laplacian operator in \cite{paquin2022asymptotics}, which relied on an explicitly known formula for the source-neutral Green's function as well as a rapidly converging sum representation of the Helmoholtz Green's function. For the latter, modifications to \eqref{Gdecompose} would need to be made to the construction of the Green's function; depending on the value of $\alpha$, there may be additional, weaker terms other than $u_0(\bx-\bx_0)$ in \eqref{Gdecompose} that become unbounded at $\bx_0$. These terms must be included in the original decomposition so that the remaining regular part to be computed numerically remains bounded at $\bx_0$. To do this, one would consider the free space Helmholtz Green's function, $G_f(\bx)$ satisfying
\BE \label{Gfreehelm}
	-(-\Delta)^\alpha G_f - v G_f = \delta(\bx) \,, \qquad \bx \in \mathbb{R}^2 \,,
\EE
for some constant scalar $v$. The Fourier transform of $G_f(\bx)$, which we denote $\hat{G}_f(\bk)$ for the Fourier variable $\bk$, satisfies
\BE \label{Gfreefourier}
	 \hat{G}_f(\bk) = -\frac{1}{v + |\bk|^{2\alpha}} \,.
\EE
For large $|\bk|$, \eqref{Gfreefourier} has the expansion
\BE \label{Gfreefourierexpand}
\hat{G}_f = -\frac{1}{|\bk|^{2\alpha}} \left\lbrack 1 - \frac{v}{|\bk|^{2\alpha}} +  \frac{v^2}{|\bk|^{4\alpha}}   - \frac{v^3}{|\bk|^{6\alpha}} + \ldots  \right\rbrack \,, \qquad |\bk| \gg 1 \,.
\EE
The singular behavior near the origin of $G_f$ can then be obtained by inverting (in the sense of distributions) the terms in the sum \eqref{Gfreefourierexpand}. The $n=1$ term gives rise to $u_0$ in \eqref{Gdecompose} with $c_\alpha$ defined in \eqref{ubar}, while the $n>1$ terms are successively weaker singular terms whose boundedness at $\bx = \bzero$ depends on the value of $\alpha$. The computation of the Helmholtz Green's function would also be useful in determining the full distribution of search times where the constant $v$ in \eqref{Gfreehelm} plays the role of the Laplace transform variable -- we refer to, e.g., \cite{lindsay2016hybrid, bressloff2021asymptotic, grebenkov2019full} for more details.

We now comment on limitations of our work. As alluded to above, our method of images-type approach to computing the Neumann Green's function of $\mAn$ was made possible by the square domain. While only slight modifications are required for rectangular domains, our approach does not extend naturally to general bounded domains. In such settings, it is also unclear, based on the argument presented in \S \ref{sec:reflect}, what particle-boundary interaction is being modeled when homogeneous Neumann boundary conditions are imposed via the spectral interpretation of the fractional Laplacian. On general domains, adapting our approach to the Neumann-type boundary conditions discussed in \cite{parks2025nonlocal, dipierro2017nonlocal} may be an interesting and worthwhile endeavor.

A second limitation in the computation of $G_\alpha^{(n)}(\bx;\bx_0)$ arises when $\bx_0$ is located near the reflecting boundary. Since the support of the cutoff function $\chi$ must lie entirely in the domain, the method may require a finer grid to resolve the steep gradients that may result. A third limitation is that our asymptotic analysis relies on a solution of the inner problem \eqref{Ui0eq}, or at least its far-field behavior, to the first two leading terms. While our assumption of the target as a circular disk allowed for an explicit solution, a more general method is needed in order to handle general target geometries. For the analogous problem with Brownian motion in two dimensions, conformal mapping \cite{ransford1995potential} and boundary integral \cite{dijkstra2008numerical} methods have been employed for non-circular geometries. It may be worthwhile to develop an analogous framework for the problem involving the fractional Laplacian.

\bibliographystyle{siam}
\bibliography{fraclap_regular}

\end{document}